\documentclass[a4paper]{amsart}
\usepackage[margin=3cm]{geometry}
\usepackage[utf8]{inputenc}
\usepackage[english]{babel}
\usepackage{amssymb}
\usepackage{amsmath}
\usepackage{amsthm}
\usepackage{amsfonts}
\usepackage{array}
\usepackage{comment,stmaryrd}
\usepackage{dirtytalk}
\usepackage{enumitem}
\usepackage{hyperref}
\usepackage{mathrsfs}
\usepackage{tikz}
\usetikzlibrary{arrows,decorations.markings}
\usepackage{tikz-cd, mathtools}
\usepackage{mathabx}
\usepackage{xfrac}
\usepackage[aligntableaux=center]{ytableau}

\usepackage{graphicx}
\graphicspath{ {./Images/} }

\addto\captionsenglish{}

\theoremstyle{plain}
\newtheorem{theorem}{Theorem}[section]
\newtheorem{deftheorem}[theorem]{Definition-Theorem}
\newtheorem{lemma}[theorem]{Lemma}

\newtheorem{proposition}[theorem]{Proposition}
\newtheorem{corollary}[theorem]{Corollary}

\theoremstyle{definition}
\newtheorem{definition}[theorem]{Definition}
\newtheorem{example}[theorem]{Example}

\theoremstyle{remark}
\newtheorem{remark}[theorem]{Remark}

\newtheorem*{theorem*}{Theorem}
\newtheorem*{conjecture*}{Conjecture}

\DeclareMathOperator{\End}{End}

\DeclareMathOperator{\Sym}{Sym}

\DeclareMathOperator{\Spec}{Spec}

\DeclareMathOperator{\Blow}{Bl}

\DeclareMathOperator{\Rep}{Rep}

\DeclareMathOperator{\Hilb}{Hilb}
\DeclareMathOperator{\Quot}{Quot}

\newcommand\Part[1][]{\mathcal{P}_{#1}} %partitions
\newcommand\CorePart[1]{\mathcal{P}_{#1, \text{\rm core}}} %core partitions
\newcommand\PartRect[1]{\mathcal{R}_{#1}} %partitions inside a rectangle

\newcommand\DPart[1][]{\mathcal{W}_{#1}} %Type D Young walls
\newcommand\DCorePart[1]{\mathcal{W}_{#1, \text{core}}} %Type D core Young walls
\newcommand\DPartFull[1][]{\mathcal{W}^{\text{f}}_{#1}} %Type D full Young walls
\newcommand\DPartHalf[1][]{\mathcal{W}^{\text{h}}_{#1}} %Type D half Young walls

\DeclareMathOperator{\Aff}{\mathbb{A}}
\DeclareMathOperator{\Prj}{\mathbb{P}}
\DeclareMathOperator{\GL}{GL}

\DeclareMathOperator{\PSL}{PSL}
\DeclareMathOperator{\SL}{SL}
\DeclareMathOperator{\Sp}{Sp}
\DeclareMathOperator{\Z}{\mathbb{Z}}

\DeclareMathOperator{\N}{\mathbb{N}}

\DeclareMathOperator{\C}{\mathbb{C}}

\DeclareMathOperator{\SO}{\mathrm{SO}}
\DeclareMathOperator{\SU}{\mathrm{SU}}

\newcommand{\CP}{\mathcal{P}}		% set of partitions
\newcommand{\wt}{\mathrm{wt}}		% weight map
\newcommand{\BF}{\mathbf{F}}                    % better Fock space

\newcommand{\ideal}{\mathfrak{a}}

\newcommand{\vect}[1]{|#1\rangle}

\makeatletter
\let\@wraptoccontribs\wraptoccontribs
\makeatother

\title[Kleinian singularities]{Kleinian singularities: some geometry, \\ combinatorics and representation theory}

\author{Lukas Bertsch}
\address{University of Vienna, Austria}
\email{lukas.bertsch@univie.ac.at }

\author{\'Ad\'am Gyenge}
\address{Budapest University of Technology and Economics, Department of Algebra and Geometry, M\H{u}egyetem rakpart 3-9., 1111, Budapest, Hungary}
\email{Gyenge.Adam@ttk.bme.hu}
\author{Bal\'{a}zs Szendr\H{o}i}
\address{University of Vienna, Austria}
\email{balazs.szendroi@univie.ac.at}

\begin{document}

\begin{abstract} We review the relationship between discrete groups of symmetries of Euclidean three-space, constructions in algebraic geometry around Kleinian
singularities including versions of Hilbert and Quot schemes, and their relationship to finite-dimensional and affine Lie algebras via the McKay correspondence. 
We focus on combinatorial aspects, such as the enumeration of certain types of partition-like objects, reviewing in particular a recently developed root-of-unity-substitution calculus. While the most complete results are in type $A$, we also develop aspects of the theory in type $D$, and end with some questions about the exceptional type $E$ cases. 
\end{abstract}

\dedicatory{Dedicated to Felix Klein on the 100th anniversary of his death}

\maketitle

\section*{Introduction}\label{sec:intro}

In the Preface to his 1884 {\it Vorlesungen \"uber das Ikosaeder}~\cite{KleinD}, Felix Klein writes 

\begin{quote}
{\it [Ich m\"ochte] dem Herr Prof.~Lie in Christiana [...] meinen besonderen Dank aussprechen. Meine Verpflichtungen gegen Hrn.~Lie gehen in die Jahre 1869-70 zur\"uck, wo wir in engem Verkehre mit einander unsere Studienzeit in Berlin und Paris abschlossen. Wir fassten damals gemeinsam den Gedanken, \"uberhaupt solche geometrische oder analytische Gebilde in Betracht zu ziehen, welche durch Gruppen von Aerderungen in sich selbst transformirt werden. Dieser Gedanke ist f\"ur unsere beiderseitigen sp\"ateren Arbeiten, soweit dieselben auch auseinauder zu liegen scheinen, bestimmend geblieben. W\"ahrend ich selbst in erster Linie Gruppen discreter Operationen ins Auge fasste und also insbesondere zur Untersuchung der regulieren K\"orper und ihrer Beziehung zur Gleichungstheorie gef\"uhrt wurde, hat Hr.~Lie von vorneherein die schwierigere Theorie der continuirlichen Transformationsgruppen [...] in Angriff genommen.}\footnote{{\it I [wish to] express my special thanks to my honoured friend Professor Lie in Christiania [...] My indebtedness to Professor Lie dates back to the years 1869-70, when we were spending the last period of our student-life in Berlin and Paris together in intimate comradeship. At that time we jointly conceived the scheme of investigating geometric or analytic forms susceptible of transformation by means of groups of changes. This purpose has been of direct influence in our subsequent labours, though these may have appeared to lie far asunder. Whilst I primarily directed my attention to groups of discrete operations, and was thus led to the investigation of regular solids and their relations to the theory of equations, Professor Lie attacked the more recondite theory of continued groups of transformations...}~\cite{Klein}.} 
\end{quote}

One cannot but be in awe of those ``intimate'' conversations, out of which grew large parts of modern mathematics, including Klein's Erlangen Programme as well as the vast subjects of Lie groups and Lie algebras. While indeed the study of discrete and that of continuous groups of symmetries may have appeared to ``lie far asunder'' for a while, later investigations have firmly brought these two fields close together again.

Our aim in this review is to present one aspect of this relationship, a theme that connects discrete groups of symmetries of Euclidean three-space to algebraic geometry around so-called Kleinian singularities, and their relationship to finite-dimensional and affine Lie algebras via the McKay correspondence. Our emphasis will be more strongly combinatorial than earlier reviews~\cite{ItoNakamura, ReidLacorrespondance, CrawReview}, spending some time on the combinatorics of partitions and related objects and their enumeration, including a discussion of some recent work~\cite{gyenge2018euler, BGS} on a curious root-of-unity-substitution calculus for enumerating labelled partition-like objects. 

One topic we will not discuss in this review is the close relationship of the subject to Nakajima quiver varieties~\cite{KronheimerNakajima, Nakajima, NakajimaBook, MaulikOkounkov}. Many of the algebraic varieties that we are going to meet, in particular the Hilbert and Quot schemes, can be realised as Nakajima quiver varieties~\cite{Kuznetsov,CGGSHilbPaper, CGGSQuotPaper, CrawYamagishi, Paegelow}. However, introducing all the necessary notation and technology would have taken us too far afield. The recent review by Craw~\cite{CrawReview} complements ours well in this regard. 

We will only touch on some aspects of the relevant representation theory: apart from representations of finite groups, only those of affine Lie algebras will make an appearance. More exotic constructions such as Cherednik algebras and Yangians will not be covered; see for example \cite{GordonICM, MaulikOkounkov} and references therein. 

In Section~\ref{sec:basic}, we introduce the basic actors of our story: the finite subgroup $G<\SL(2,\C)$, the associated Kleinian singularity $X=\Aff^2/G$ and its resolution, aspects of the equivariant geometry, and finally Hilbert and Quot schemes. We have chosen to treat the well known, early part of the theory in a little more detail than perhaps necessary, in order to emphasise Klein's contributions. In Section~\ref{sec:abelian}, which is really the heart of our review, the case when the group $G$ is abelian is treated in detail, with a strong emphasis on the associated combinatorics. In Section~\ref{sec:nonabelian}, we discuss how the story generalises to nonabelian~$G$. 

Our key actor will be the finite subgroup $G<\SL(2,\C)$. The questions studied here could be generalised in (at least) three different directions: one can consider 
\begin{enumerate}
\item a finite subgroup $G<\GL(2,\C)$, the 2-dimensional non-Calabi--Yau case; 
\item a finite subgroup $G<\SL(n,\C)$ for $n>2$, the higher dimensional Calabi--Yau case; and
\item a finite subgroup $G<\Sp(2n,\C)$ for $n>1$, the higher dimensional holomorphic symplectic case. 
\end{enumerate}
We will briefly comment on some examples of (1)  in~\ref{sec:GL2}. We will not discuss (2)-(3) at all; see for example~\cite{ItoNakajima, FuSurvey, DavisonOngaroSzendroi}.

\subsection*{Notation} We work over the complex number field. For a non-trivial finite group~$G$, we will always use the convention that the number of conjugacy classes of $G$ will be $1+r$ for some integer $r>0$. Let $\Rep(G)$ denote the set of all finite-dimensional representations of~$G$ over~$\C$ up to isomorphism, and let $\rho_0,\ldots, \rho_r\in\Rep(G)$ be a complete list of distinct irreducibles, with underlying vector spaces $V_0, \ldots, V_r$, $\rho_0$ being the trivial representation. We will denote $I=\{0,\ldots, r\}$ the index set of irreducibles. The letter $V$ will usually denote the vector space underlying some fixed representation $\rho\in\Rep(G)$ of $G$; for $G<\SL(2,\C)$, this will always be the ``given'' two-dimensional representation of~$G$ coming from this embedding.

\subsection*{Acknowledgements} We would like to thank Alastair Craw, S\o ren Gammelgaard, Rapha\" el Paegelow and Michael Schlosser for comments. Á.Gy.~was supported by a János Bolyai Research Scholarship of the Hungarian Academy of Sciences and by the  “Élvonal (Frontier)” Grant KKP 144148. 

\section{The basic characters of our story}\label{sec:basic}

\subsection{Some groups of rotations}\label{sec:Subgroups}

Our point of departure is the group $\SO(3)$ of all orientation-preserving isometries of three-dimensional Euclidean space with a fixed point at the origin. As is well known, each element of $\SO(3)$ is a 
rotation around an axis, a line in three-dimensional space passing through the origin. Finite subgroups of $\SO(3)$ arise naturally as rotational symmetry groups of geometric objects. Klein~\cite{Klein} described the full set of finite subgroups of $\SO(3)$ up to conjugation as follows.

\begin{enumerate}
\item Cyclic group $C_m=\langle s \mid s^m = e \rangle <\SO(3)$: the group of rotational symmetries of the perpendicular pyramid over a regular planar $m$-gon. 
\item Dihedral group $D_m=\langle s, t \mid s^2 = t^2 = (st)^m = e \rangle <\SO(3)$: the group of rotational symmetries of the perpendicular double pyramid over a regular planar $m$-gon. 
\item Tetrahedral group $\Gamma_T<\SO(3)$: the group of rotational symmetries of the regular tetrahedron. It can be presented as
\[ \Gamma_T \cong \langle s, t \mid s^2 = t^3 = (st)^3 = e \rangle \; . \]
Here $t$ can be taken to be a rotation around an axis containing a vertex and the midpoint of the opposite face, and then $s$ is a rotation around an axis containing midpoints of opposite edges. 
As an abstract group, the group $\Gamma_T$ is isomorphic to the alternating group $A_4$, as can be seen by considering its action on the vertices of the tetrahedron. 
\item Octahedral group $\Gamma_O<\SO(3)$: the group of rotational symmetries of the cube and the (regular) octahedron. Recall that the cube and the octahedron are duals, meaning that midpoints of faces of one regular solid give vertices of the other, and hence their symmetry groups are the same. This group can be presented as
\[\Gamma_O \cong \langle s, t \mid s^2 = t^3 = (st)^4 = e \rangle \; . \]
As an abstract group, $\Gamma_O$ is isomorphic to the symmetric group $S_4$, as can be seen by considering its action on the big diagonals of the cube. 
\item Icosahedral group $\Gamma_I<\SO(3)$: the group of rotational symmetries of the (regular) icosahedron and the (regular) dodecahedron, consisting of $60$ elements. The icosahedron and the dodecahedron are also duals. The group $\Gamma_I$ can be presented as
\[ \Gamma_I \cong \langle s, t \mid s^2 = t^3 = (st)^5 = e \rangle \; . \]
As an abstract group, $\Gamma_I$ is isomorphic to the alternating group $A_5$, as can be seen by considering its action on the set of vertex-embedded regular tetrahedra of the dodecahedron. 
\end{enumerate}

As discussed by Klein, one can realize finite subgroups of $\SO(3)$, such as the symmetry groups of the Platonic solids, as finite subgroups of $\PSL(2, \mathbb{C})$ using Möbius transformations acting on the Riemann sphere $\C\cup\{\infty\}$. 
%Namely, an element \[\left[\begin{pmatrix} a & b \\ c & d \end{pmatrix}\right]\in \PSL(2, \mathbb{C})\] acts on a point $z\in\C\cup\{\infty\}$ in the Riemann sphere by  $f(z) = \frac{az + b}{cz + d}$.
We get a diagram
$$\begin{tikzcd} \SU(2) \arrow[r,hook] \arrow[d, "2:1"] & \SL(2,\C) \arrow[d, "2:1"] \\ \SO(3) \arrow[r,hook] & \PSL(2,\C)
\end{tikzcd}$$

Given a finite subgroup $\Gamma<\SO(3)$, pulling back by the double cover map one gets a finite subgroup $\widetilde \Gamma<\SU(2)$, the corresponding ``binary'' group. The order of this group is $|\widetilde\Gamma|=2|\Gamma|$. This construction gives essentially all finite subgroups of $\SU(2)$; there is a small discrepancy in the abelian case, where this construction only gives the abelian (in particular cyclic) subgroups of $\SU(2)$ of even order.  

On the other hand, given a finite subgroup $G<\SL(2,\C)$, one can take an arbitrary hermitian form on $\C^2$, and then average it out over $G$, to get $G<\SU(2)$ in an appropriate basis. Putting all these facts together, we get the following result. 
\begin{theorem}[Klein~\cite{Klein}] \label{thm:classification} The following is a complete list of finite subgroups of $\SL(2,\C)$ up to conjugation.
\begin{enumerate}
\item The cyclic subgroup generated by
\[ \sigma = \begin{pmatrix}
\omega & 0 \\
0 & \omega^{-1}
\end{pmatrix}, \]
where $ \omega = e^{2\pi i / m} $ is a primitive $m$-th root of unity.
For odd $m$, this group maps isomorphically onto $C_m<\SO(3)$. For even $m=2n$, it provides the ''binary cyclic group'' $\widetilde C_n<\SU(2)$ with a double cover map to $C_n<\SO(3)$.
\item The Binary Dihedral Group $\widetilde D_m$, generated by
\[ \sigma = \begin{pmatrix}
\omega & 0 \\
0 & \omega^{-1}
\end{pmatrix},\quad
\tau = \begin{pmatrix}
0 & 1 \\
-1 & 0
\end{pmatrix}, \]
where $ \omega = e^{2\pi i / m} $ is a primitive $m$-th root of unity. \item The Binary Tetrahedral Group $\widetilde \Gamma_T$, generated by
\[
\sigma = \begin{pmatrix}
i & 0 \\
0 & -i
\end{pmatrix}, \quad
\tau = \frac{1}{\sqrt{2}}\begin{pmatrix}
1 & 1 \\
-1 & 1
\end{pmatrix}.
\]
\item The Binary Octahedral Group $\widetilde \Gamma_O$, generated by
\[
\sigma = \frac{1}{\sqrt{2}}\begin{pmatrix}
1 & i \\
i & 1
\end{pmatrix}, \quad
\tau = \frac{1}{\sqrt{2}}\begin{pmatrix}
1 & 1 \\
-1 & 1
\end{pmatrix}.
\]

\item The Binary Icosahedral Group $\widetilde  \Gamma_I$, generated by
\[
\sigma = \begin{pmatrix}
\phi & 1 \\
1 & -\phi
\end{pmatrix}, \quad
\tau = \begin{pmatrix}
\omega & 0 \\
0 & \omega^{-1}
\end{pmatrix},
\]
where $ \phi = \frac{1 + \sqrt{5}}{2} $ is the golden ratio, and $ \omega = e^{2\pi i / 5} $ is a primitive fifth root of unity.
\end{enumerate}
\end{theorem}

For the rest of this review, we will focus on these binary groups, so we will henceforth reserve the notation $G$
to denote an arbitrary finite subgroup $G< \SL(2,\C)$, which can be abelian, (binary) dihedral, or exceptional (corresponding to the symmetry group of one of the regular solids). 

\subsection{Invariant theory and algebraic geometry}\label{subsec:singularity}

Consider the complex affine plane~$\Aff^2$, with coordinate algebra $\C[\Aff^2] \cong \C[x,y]$. A finite subgroup $G < \SL(2,\C)$ acts on the two-dimensional vector space $V=\langle x,y\rangle$, which extends to an action on $\C[x,y]=\Sym^*V$ by algebra automorphisms. Dually, this gives a geometric action on~$\Aff^2=\Spec\Sym^*V$. 

Given the action of $G$ on $\C[\Aff^2]$, it was very natural for Klein to consider, in the spirit of 19th century invariant theory, the invariant algebra $\C[\Aff^2]^{G}$.
\begin{theorem}[Klein~\cite{Klein}] \label{thm gens} Given a finite subgroup $G < \SL(2,\C)$, the invariant algebra $\C[\Aff^2]^{G}$ is always generated by three invariant polynomials $u, v, w\in \C[\Aff^2]^{G}$, satisfying a single polynomial relation $g(u, v, w)=0$, with the lowest order term in $g$ being quadratic.  
\end{theorem}
We obtain the algebraic description 
\[
\C[\Aff^2]^{G} \cong \C[u, v, w]/\langle g \rangle
\]
of our invariant algebra. The corresponding {\em Kleinian singularity} $X = \Aff^2 / G=\Spec \C[\Aff^2]^{G}$ is the affine variety whose coordinate algebra is the invariant algebra $\C[X] = \C[\Aff^2]^{G}$; the quotient map $\Aff^2 \to X$ corresponds to the inclusion of algebras $\C[\Aff^2]^{G} \subset \C[\Aff^2]$. The variety $X$ has a unique singular point at the image
of the origin $0=[(0,0)]\in X$. Theorem~\ref{thm gens} implies also that a Kleinian singularity $X$ can always be embedded in affine three-space as a hypersurface. Kleinian singularities are also known as rational double points or du Val singularities, for various reasons too long to explain here; they are ubiquitous in algebraic geometry and singularity theory~\cite{Durfee}. 

\begin{figure}

\begin{center}
    \includegraphics[height=5cm]{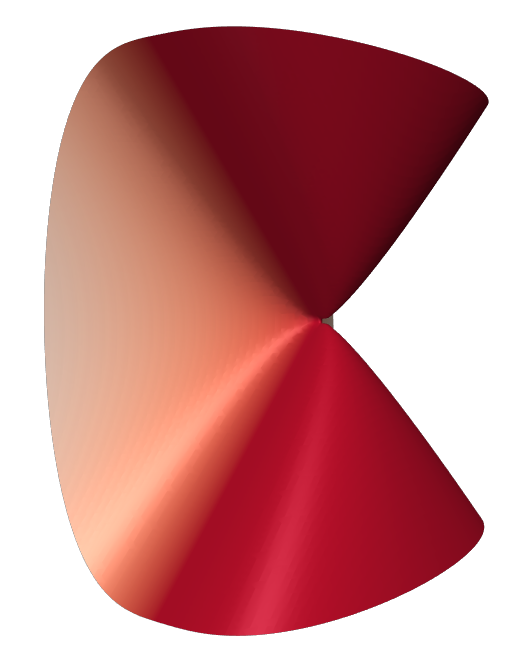}
    \caption{The real locus of the singular surface $u^2 - w^2 - v^3=0$}\label{fig:A2sing}
\end{center}
\end{figure}

\begin{example}\label{ex:A2-invariant-ring} Consider $G \cong C_3=\langle\sigma\rangle$ acting on~$V=\langle x,y\rangle$ by $$\sigma\circ \left(\begin{array}{c} x \\ y \end{array}\right) = \left(\begin{array}{cc} \zeta & 0  \\0 & \zeta^{-1}\end{array}\right) \left(\begin{array}{c} x \\ y \end{array}\right), $$ where $\zeta$ is a primitive $3$-rd root of unity. In this case the algebra $\C[\Aff^2] = \C[x,y]$ has a basis of monomials $x^a y^b$, all of which are eigenvectors of this action. 
In particular, the invariant ring is spanned by the monomials $x^a y^b$ with $a - b \equiv 0 \pmod{3}$. These monomials form a monoid under multiplication, which is generated by $u = x^3$, $v = xy$ and $w = y^3$, satisfying the relation $g(u,v,w) = u w - v^3=0$. After convincing ourselves that $g$~generates all relations among these generators, we obtain a description of the invariant ring as $$\C[\Aff^2]^{G} \cong \C[u,v,w]/(u w - v^3) \; .$$ After a simple change of coordinates, one obtains the  alternative form $$\C[\Aff^2]^{G} \cong \C[u,v,w]/(u^2 - w^2 - v^3) \; .$$
See Figure~\ref{fig:A2sing} for the real point set of this singularity. 

For a detailed study of the smallest dihedral example, see~\cite[1.3-1.4]{ReidADE}.
\end{example}

\vspace{.1in}

\setlength{\extrarowheight}{3pt}
\begin{center}
\begin{tabular}{c|c|c} Finite subgroup $G<\SL(2,\C)$ & Polynomial defining Kleinian singularity & Dynkin type \\[3pt] \hline Cyclic $C_{r+1} \ (r\geq 1)$& $x^2 + y^2 + z^{r+1}$ & $A_r$ \\
Binary dihedral $\widetilde D_{r-2} \ (r\geq 4)$  & $x^2+y^2z + z^{r-1}$ & $D_r$\\ 
Binary tetrahedral  $\widetilde \Gamma_{T}$ & $x^2+y^3+z^4$ & $E_6$ \\
Binary octahedral $\widetilde \Gamma_{O}$ & $x^2+y^3+yz^3$ &   $E_7$\\
Binary icosahedral $\widetilde \Gamma_{I}$ & $x^2+y^3+z^5$ &   $E_8$\end{tabular}
\end{center}

\vspace{.1in}

In order to study singular quotient varieties such as Kleinian singularities $X=\Aff^2/G$, modern algebraic geometry offers two different tools, both of which allow us to replace the singular variety~$X$ by a smooth space.
\begin{enumerate}
    \item A \textit{resolution of the singularity on} $X$ is a smooth quasiprojective variety $Y$ together with a proper map $\pi\colon Y \to X$ which is an isomorphism over the smooth locus of $X$.
    \item Instead of forming the quotient in the first place, one can consider the \textit{equivariant geometry} of $\Aff^2$ with respect to its action of $G$. This is equivalent to studying the quotient orbifold or stack $[\Aff^2/G]$.
\end{enumerate}

A central theme in the geometry of Kleinian singularities is that there is a strong relationship between the two viewpoints. A specific form of this relationship was discovered by John McKay~\cite{McKay} in the form of an (almost) matching between two Dynkin diagrams that arise naturally in the two situations. It is therefore known as the \textit{McKay correspondence}. For the rest of this section, we briefly discuss the resolution of a Kleinian singularity. In the next section, we will elaborate on the equivariant geometry, and recall McKay's observation.

\begin{theorem}\label{thm:kleinian-resolution} Given a finite subgroup $G<\SL(2,\C)$, the Kleinian singularity $X=\Aff^2/G$ has a unique minimal resolution $\pi\colon Y\to X$, obtained by iteratively blowing up the singular locus. The resolution $\pi$ has an exceptional divisor $\pi^{-1}(0)\subset Y$ that is a union of rational curves intersecting transversally, with dual graph a Dynkin diagram. This diagram is of type A, D or E respectively for cyclic groups $C_m$, binary dihedral groups $\widetilde D_m$, and binary Platonic groups $\widetilde T, \widetilde O, \widetilde I$.  
\end{theorem}

\begin{example}\label{ex:A2-resolution}
    From Example \ref{ex:A2-invariant-ring}, we know that the singularity for $G\cong C_3$ has coordinate ring $$\C[\Aff^2]^{G} \cong \C[u, v, w] / (u w - v^3) \; .$$ In this case, a short calculation shows that the blowup $\Blow_0(X)$ of $X = \Spec \C[\Aff^2]^{G}$ at the origin is already smooth, and therefore gives a resolution of singularities $$\pi\colon Y=\Blow_0(X) \rightarrow X \; .$$ The exceptional divisor of~$\pi$ is the projectivised normal cone of $0\in X$, which is given by the lowest-degree part of the equation as $\{u w = 0\}\subset \Prj^2$. This defines the union of two projective lines which intersect transversally in one point, giving us a dual graph with two vertices connected by an edge, the $A_2$ Dynkin diagram.  
\end{example}

The type of the Dynkin diagram from Theorem~\ref{thm:kleinian-resolution} is also called the type of the subgroup $G<\SL(2,\C)$. We thus see that the subgroup $G\cong C_3$ of $\SL(2,\C)$ is of type $A_2$. More generally, $C_{r+1} \cong G<\SL(2,\C)$ corresponds to type $A_r$. 

\subsection{Equivariant geometry and the McKay Quiver}\label{sec:Equivariant}

The algebraic study of the geometry of an algebraic variety involves studying (coherent) sheaves over the variety. If we are interested in the equivariant geometry of a variety $Z$ with respect to some group $G$ acting on it, we need the notion of an {\it equivariant sheaf} on $(Z,G)$, incorporating the action. Here we give the definition in the case of an affine variety $Z$ with coordinate ring $R=\C[Z]$, where the study of sheaves on $Z$ reduces to that of modules over the coordinate ring $R$. 

\begin{definition}
    Consider a commutative $\C$-algebra $R$ and a finite group $G$ acting on $R$, written $G \times R \to R, \; (g,f) \mapsto g(f)$. An \textit{equivariant (left) module} over $R$ with respect to this action is a $R$-module $M$ together with a $\C$-linear (left) action of $G$ on $M$ such that \begin{equation}\label{eq:equiv-module}g(fm) = g(f)g(m) \quad \text{for all} \; g \in G, \; f \in R, \; m \in M \; .\end{equation}
\end{definition}

Hence the structure of an equivariant module on a vector space $M$ is given by two intertwined module structures: one over $R$, and one over the group algebra $\C G$. One can define a new algebra incorporating both of these actions.

\begin{definition}
    The \textit{skew-group algebra} of $R$ with respect to the action of $G$ is the associative algebra $R \rtimes G$ with underlying vector space $R \otimes \C G$ and multiplication $$(f \otimes g)(f' \otimes g') = (f g(f')) \otimes (g g') \; .$$
\end{definition}

One can check that this multiplication law is indeed associative. Furthermore, both $R$ and $\C G$ are subalgebras of $R \rtimes G$, embedded as $f \mapsto f \otimes 1$ and $g \mapsto 1 \otimes g$. Under this embedding we have $fg = (f \otimes 1) (1 \otimes g) = f \otimes g$, so we can omit the symbol $\otimes$ from notation. The relation \begin{equation}\label{eq:skew-commutator} gf = g(f) g \end{equation} holds in $R \rtimes G$, which shows that $R \rtimes G$ is non-commutative, even if both $R$ and $G$ are commutative. The commutation relation (\ref{eq:skew-commutator}) reflects the condition (\ref{eq:equiv-module}) in the action on a module and therefore ensures the following.

\begin{proposition} 
Given a vector space $M$, a $G$-equivariant $R$-module structure on $M$ is equivalent to that of a left $R \rtimes G$-module structure on $M$.
\end{proposition}

Hence, given an affine variety $Z$ with coordinate ring $R=\C[Z]$ together with an action of a finite group $G$, the equivariant geometry of $(Z,G)$ is equivalent to the \textit{non-commutative geometry} of the algebra $R \rtimes G$. Note that $R$, which is itself naturally a $G$-equivariant $R$-module, can be identified as $$R \cong (R \rtimes G)e_0 \; ,$$ where $e_0 = \frac{1}{|G|}\sum_{g \in G}g \in \C G$ is the invariant idempotent. The invariant ring can be recovered as the (non-unital) subring $$R^{G} \cong e_0(R \rtimes G)e_0 \; .$$ There is also a second way in which $R^{G}$ is a subring, given by the well-known
\begin{proposition}
    Suppose that $R$ is an integral domain and the action of $G$ on $R$ is faithful. Then $$Z(R \rtimes G) = R^{G} \; .$$
\end{proposition}
\noindent This inclusion of $R^{G} = Z(R \rtimes G)$ as a subring in $R \rtimes G$ is different from the inclusion as $R^{G} \cong e_0 (R \rtimes G) e_0$ in $R \rtimes G$.

The skew-group algebra $R \rtimes G$ thus arises naturally from the viewpoint of equivariant affine geometry. One can, however, think about it the other way around also: consider first the group~$G$, whose representation theory is well-understood. Then take generators of $R$, and consider their action on representations of $G$. The starting point for this perspective is the following standard result.

\begin{theorem}\cite[Prop.\ 3.29]{fultonharris1991} \label{thm:group-algebra-structure}
    Let $G$ be a finite group, $V_0, \ldots, V_r$ vector spaces underlying its different irreducible representations $\rho_0, \ldots, \rho_r$.
    Then there is an isomorphism of $\C$-algebras
     $$\C G \cong \prod_{i = 0}^r \End_{\C}(V_i) \; .$$ 
\end{theorem}

It follows from this result that the algebra $\C G$ is Morita equivalent to (has equivalent category of modules with) the commutative algebra $\prod_{i = 0}^r \C$. If $G$ is abelian, we even have a ring isomorphism $\C G \cong \prod_{i = 0}^r \C$. In geometric language, the representation theory of $G$ is in a certain sense equivalent to that of the coordinate ring of $(r+1)$ points. To this picture we now add the data of how $R$ acts on representations of $G$.

Suppose that $R$, as an algebra with $G$-action, is generated by a finite-dimensional representation $V$ of $G$. Then a $G$-equivariant action of $R$ on a $G$-representation $M$ is specified by the action of $V$ on $M$, which is by (\ref{eq:equiv-module}) a $G$-intertwiner $$V \otimes M \to M \; .$$ Suppose $m \in M$ is a vector that lies inside an irreducible subrepresentation isomorphic to $V_i$, and $f \in V$, then $f m$, being the image of $f \otimes m$ under this multiplication map, lies in the quotient of a representation isomorphic to $V \otimes V_i$. In geometric language, if $m$ is supported on the point labeled $i$, then $f m$ is supported on the points corresponding to irreducible representations appearing in $V \otimes V_i$. This observation motivates the following definition; recall that a {\it quiver} is an oriented graph, consisting of a set of vertices and oriented edges (arrows) between them. 

\begin{definition}[McKay~\cite{McKay}]
    Let $G$ be a finite group with irreducible representations $\rho_0,\ldots, \rho_r$ with underlying vector spaces $V_0, \ldots, V_r$, and a fixed representation~$\rho$ on a vector space~$V$. The \textit{McKay quiver} $Q$ of the pair $(G, \rho)$ is defined as follows.
    \begin{itemize}
        \item $Q$ has vertex set $I = \{0, \ldots, r\}$, in one-to-one correspondence with the set of irreducible representations $\{\rho_0,\ldots, \rho_r\}$ of~$G$.
        \item For vertices $i,j \in I$, the number of oriented edges from $i$ to $j$ equals the multiplicity of the $G$-representation $V_j$ in $V \otimes V_i$.
    \end{itemize}
\end{definition}

We return to our setting where $G < \SL(2,\C)$ is a finite subgroup and $\rho$ is the given two-dimensional representation on the vector space $V=\langle x, y\rangle$. The multiplicity of $V_j$ in $V \otimes V_i$ in this case is the same as that of $V_i$ in $V^* \otimes V_j \cong V \otimes V_j$, so the arrows of the McKay quiver of $(G,V)$ come in opposite pairs in this case.

The McKay quiver $Q$ encodes information about both the group $G$ and its action on $V$. Hence, when we consider the commutative algebra $R$ freely generated by $V$, we would expect to be able to reconstruct the algebra $R \rtimes G$ from the quiver $Q$. This is indeed the case. Consider the path algebra $\C Q$ of $Q$: the set of paths in $Q$ form a basis of $\C Q$, and multiplication is given by concatenation of paths (when the endpoint of the first path does not match the starting point of the second, their product is zero). Notice that $\C Q$ contains $\prod_{i = 0}^r \C_i$ as the subalgebra of paths of length zero.

\begin{theorem}[Reiten and van den Bergh~\cite{RvdB}]
%[Prop.2.4,Thm.5.6]
\label{thm:morita-equivalence}
    There is a two-sided ideal $\mathfrak{p}$ in $\C Q$ such that the quotient $$\Pi = \C Q / \mathfrak{p}$$ is Morita equivalent to $R \rtimes G$. When $G$ is abelian, there is even an isomorphism $\Pi \cong R \rtimes G$.
\end{theorem}

The quotient $\Pi$ is called the \textit{preprojective algebra}; for a detailed definition, see e.g.\ \cite[Def.\ 5.2]{kirillov2016quivers} or~\cite[Section 3.2]{CrawReview}. In essence, the fact that $R$ is commutative imposes relations in $\C Q$ between paths of length two, and $\mathfrak{p}$ is generated by these relations. Notice that, since the center of an algebra is a Morita invariant, we also have $R^{G} = Z(R \rtimes G) \cong Z(\Pi)$.

\begin{example}\label{ex:A2-quiver}
    Consider the group $C_3 \cong G_{A_2} < \SL(2,\C)$. The given representation $V$ decomposes into irreducibles as $V = V_1 \oplus V_2$ where $V_1 = \langle x\rangle$ and $V_2 = \langle y\rangle$. Together with the trivial reprentation~$V_0$, these form the complete set of irreducible representations of $G_{A_2}$. The McKay quiver looks as follows.
    \begin{center}\begin{tikzpicture}[scale=0.6, font=\footnotesize, fill=black!20, auto]
        \node[circle,draw] (a) at (180:2) {$\rho_0$};
        \node[circle,draw] (b) at (60:2) {$\rho_1$};
        \node[circle,draw] (c) at (300:2) {$\rho_2$};
        \draw[shorten <=2, shorten >=2, ->, bend left, bend angle=60] (a.50) to node[swap] {$x_0$} (b.190);
        \draw[shorten <=2, shorten >=2, ->, bend left, bend angle=60] (b.290) to node[swap] {$x_1$} (c.70);
        \draw[shorten <=2, shorten >=2, ->, bend left, bend angle=60] (c.170) to node[swap] {$x_2$} (a.310);
        \draw[shorten <=2, shorten >=2, ->, bend right, bend angle=60] (a.280) to node[swap] {$y_2$} (c.200);
        \draw[shorten <=2, shorten >=2, ->, bend right, bend angle=60] (c.40) to node[swap] {$y_1$} (b.320);
        \draw[shorten <=2, shorten >=2, ->, bend right, bend angle=60] (b.160) to node[swap] {$y_0$} (a.80);
    \end{tikzpicture}\end{center}
    Here the labels $x_i$ and $y_j$ indicate tensoring by $V_1$ and $V_2$, respectively. The preprojective relations in this case are precisely $x_iy_i=y_{i+1}x_{i+1}$, for $i\in\Z/3$.
\end{example}

It is easy to see that the McKay quiver for the cyclic group $C_{r+1}<\SL(2,\C)$ looks similar to that of Example \ref{ex:A2-quiver}, with $(r+1)$ nodes arranged in a circle. In other words, the McKay quiver of the cyclic group $C_{r+1}<\SL(2,\C)$ is the double quiver associated to an extended (affine) Dynkin diagram of type~$\widehat A_r$. The corresponding finite Dynkin diagram, formed by components of the exceptional divisor in the minimal resolution of the corresponding singularity, was found at the end of~\ref{subsec:singularity} to be of type $A_r$. The main observation of McKay was that this relation holds for all finite subgroups~$G<\SL(2,\C).$.

\begin{theorem}[McKay~\cite{McKay}] \label{thm:kleinian-mckay} Fix a finite subgroup $G<\SL(2,\C)$ with its given two-dimensional representation~$\rho$ on the vector space $V=\langle x, y\rangle$.
    The McKay quiver~$Q$ associated to $(G, \rho)$ is the double quiver associated to the extended (affine) ADE Dynkin diagram corresponding to (non-extended) Dynkin diagram formed by components of the exceptional divisor of the Kleinian singularity $X=\Aff^2/G$ (compare~Theorem \ref{thm:kleinian-resolution}). The distinguished vertex of the extended Dynkin diagram corresponds to the trivial representation~$\rho_0$.
\end{theorem}
This remarkable observation, which relates the resolution of a Kleinian singularity $\Aff^2/G$ and the representation theory of the group $G$, is what became known as the \textit{McKay correspondence}.  It has inspired a lot of work, which we do not have the space to summarise; reviews include~\cite{ItoNakamura, ReidLacorrespondance} as well as articles in~\cite{Ito}. We mention only one of the most general possible formulations, the so-called {\em derived McKay correspondence}. This is the result~\cite{KapranovVasserot} that states an equivalence of triangulated categories between the derived category of coherent sheaves on the minimal resolution~$Y$ of a Kleinian quotient singularity $X=\Aff^2/G$, and the derived category of $G$-equivariant coherent sheaves on the plane (equivalently, the derived category of coherent sheaves on the stack $[\Aff^2/G]$):
\[
  D^b(\text{Coh}(Y)) \simeq D_G^b(\text{Coh}(\Aff^2)).
\]
The triangulated category on the right is nothing but the category $D^b(\text{Mod}(\C[x,y] \rtimes G))$. 
A very substantial generalisation of this equivalence to higher dimensions, including all quotients $X=\Aff^3/G$ for $G<\SL(3, \mathbb{C})$, was proved in~\cite{BKR}, with a lot of further work afterwards, see for example~\cite{Ito}.

\subsection{Hilbert and Quot schemes}

For any algebraic variety $Z$, a fundamental object of interest is $\Hilb(Z)$, the \textit{Hilbert scheme of points on $Z$}. This is the moduli space of finite-length subschemes of~$Z$, decomposing into disjoint components 
\[\Hilb(Z)=\bigsqcup_n \Hilb^n(Z)\]
indexed by the length $n$ of the subscheme; see for example~\cite{NakajimaBook, FGAexplained}. As we will only encounter the affine situation in this review, we give the definition only in this case. So let $Z$ be an affine variety over~$\C$, with coordinate algebra $R=\C[Z]$. Then the $n$-th Hilbert scheme of $Z$ is defined to be
\[ \Hilb^n(Z)= \{ \mbox{quotient algebras } R\twoheadrightarrow R/\ideal \colon {\rm dim}_{\C}(R/\ideal) = n \}.
\]
As written, this is only a set, but it carries the structure of a quasiprojective scheme. 

Return to our situation, where $Z=\Aff^2$ is the affine plane. Then $\Hilb^n(\Aff^2)$ is a smooth (reduced) quasiprojective variety, known to admit a concrete description in terms of linear algebraic data~\cite{NakajimaBook, CrawReview}. Further, given
a finite subgroup $G < \SL(2,\C)$, $G$ acts on $\Hilb(\Aff^2)$ via its action on the coordinate ring $R$. Consider the fixed point locus $\Hilb(\Aff^2)^G$ of this action, which decomposes as follows: 
$$\Hilb(\Aff^2)^G = \bigsqcup_{\rho \in \Rep(G)} \Hilb^{\rho}(\Aff^2) \; .$$ 
Here $\Hilb^{\rho}(\Aff^2)$ parametrizes $G$-invariant subschemes of $\Aff^2$, whose coordinate ring is of $G$-represen\-tation type $\rho$.
We call this fixed point locus the \textit{equivariant Hilbert scheme of points}~\cite{ItoNakamura, Brion} of $(\Aff^2,G)$. 
The data of a representation type $\rho$ is equivalent to that of a vector $v \in \N^I$, which determines $\rho \cong \sum_{i \in I} \rho_i^{\oplus v_i}$. It is known that each of the spaces $\Hilb^{\rho}(\Aff^2)$ is smooth and irreducible when non-empty. For the abelian case, a bijective parametrisation of components in terms of combinatorial data will be discussed below in Theorem~\ref{thm:typeAmain}(i). 

In the general language developed in~\ref{sec:Equivariant}, we can think of the equivariant Hilbert scheme of an affine variety with group action $(Z,G)$, for example for $Z=\Aff^2$ acted on by $G<\SL(2,\C)$, as follows. Subschemes of $Z = \Spec(R)$ are the same as quotient algebras of $R$, in other words, $R$-quotient modules of the fixed $R$-module $R$. Further, $G$-equivariant subschemes are the same as left $(R \rtimes G)$-quotient modules of $R = (R \rtimes G)e_0$. On the other hand, subschemes in the quotient $Z / G$ are the same as $R^{G} \cong e_0(R \rtimes G)e_0$-module quotients of $e_0(R \rtimes G)e_0$. The following definition and theorem concerning a \textit{non-commutative Quot scheme} interpolates between these two definitions. 

\begin{deftheorem}\cite{CGGSQuotPaper}
\label{defthm:quotscheme}
    Let $e_i$ be the projector onto the summand $\End_{\C}(V_i)$ in Theorem \ref{thm:group-algebra-structure}. Consider a non-empty subset $J \subseteq I$, and let $e_J \coloneqq \sum_{j \in J} e_j\in R \rtimes G$. There exists a fine moduli scheme $\Quot_{G,J}(\Aff^2)$ that parametrises finite-dimensional left $e_J(R \rtimes G)e_J$-quotient modules of $e_J(R \rtimes G)e_0$. 
\end{deftheorem}

The process of obtaining $e_J(R \rtimes G)e_J$, and its module $e_J(R \rtimes G)e_0$, is sometimes called \textit{cornering}~\cite{Karmazinetal}.
Note that, by virtue of the Morita equivalence of Theorem \ref{thm:morita-equivalence}, our non-commutative Quot schemes can also be defined in terms of the preprojective algebra~$\Pi$, which also contains a corresponding set of idempotents.

\begin{example}\label{ex:quots}
The following special cases of the definition are worth mentioning: 
\begin{enumerate} 
\item For $J=I$, $e_J$ is simply the unit of the algebra~$R \rtimes G$, and we have 
\[\Quot_{G,I}(\Aff^2) = \Hilb(\Aff^2)^G.\] 
\item For $J=\{0\}$, the definition returns 
\[\Quot_{G,\{0\}}(\Aff^2) = \Hilb(\Aff^2/G),\]
the Hilbert scheme of points of the singular affine surface $X=\Aff^2/G$, studied in~\cite{CGGSHilbPaper, CrawYamagishi}.
\end{enumerate}
\end{example}

Recall that the equivariant Hilbert scheme $\Hilb(\Aff^2)^{G}$ decomposes according to representation type over $\C G$, and the ordinary Hilbert scheme $\Hilb(\Aff^2/G)$ decomposes according to dimension of the quotient ring. This generalises: $\Quot_{G,J}(\Aff^2)$ decomposes according to representation type over $e_J \C G e_J$, which is equivalent to the data of a vector $v_J \in \N^J$. Hence we have a decomposition \begin{equation}\label{eq:quot-decomposition}\Quot_{G,J}(\Aff^2) = \bigsqcup_{v_J \in \N^J} \Quot_{G,J}^{v_J}(\Aff^2) \; .\end{equation} It is known that each $\Quot_{G,J}^{v_J}(\Aff^2)$ is a quasi-projective scheme, in particular of finite type.

Finally, we describe natural maps between the Quot schemes for varying $J$. 

\begin{proposition}\cite{CGGSQuotPaper} \label{prop:deg} Suppose $J' \subseteq J$ are nonempty subsets of $I$. Then, using the natural idempotent $e_{J'} \in e_J(R \rtimes G)e_J$, the assignment of modules $M \mapsto e_{J'} M$ induces a morphism of Quot schemes 
\begin{equation}\label{eq:quot-degeneration}
p_{J,J'}\colon  \Quot_{G,J}^{v_J}(\Aff^2) \longrightarrow \Quot_{G,J'}^{v_{J'}}(\Aff^2) 
\end{equation}
for every $v_J \in \N^J$ and $v_{J'} = v_J|_{J'} \in \N^{J'}$. These morphisms are compatible under composition along consecutive inclusions $J'' \subseteq J' \subseteq J$. 
\end{proposition}

We shall refer to the morphisms $p_{J,J'}$ as \textit{degeneration maps}. A particular example of this morphism is a natural map \begin{equation}\label{eq:quot-degenerationspec} p_{I,\{0\}}\colon \Hilb(\Aff^2)^G \longrightarrow \Hilb(X)\end{equation} from the (smooth) equivariant Hilbert scheme to the Hilbert scheme of the quotient $X=\Aff^2/G$, defined in~\cite[3.4]{Brion}. In the language of ideals, this simply maps a $G$-equivariant ideal $\ideal\lhd\C[x,y]$ to the intersection
$\ideal\cap\C[x,y]^G\lhd\C[x,y]^G$. 

\section{The abelian case: partition combinatorics}\label{sec:abelian}

In this section, we study the case (1) of Theorem~\ref{thm:classification}, the cyclic subgroups $C_{r+1}\cong G_{A_r}<\SL(2,\C)$ of type $A_r$, in more detail. In this abelian case, we will be able to apply torus localisation to the Hilbert and Quot schemes, which will allow us to understand certain aspects of these spaces through the combinatorics of labelled Young diagrams (partitions).

\subsection{Combinatorics around labelled partitions}\label{subsec:part}

We begin by introducing some standard combinatorial constructions, which we will relate to Hilbert and Quot schemes in the next section. 

%\subsubsection{The unlabelled case: partitions}

Recall that a \textit{partition} $\lambda$ of a positive integer $n$ is a decomposition $n=\lambda_1+\ldots + \lambda_k$ into positive integers in weakly decreasing order. We write $\lambda=(\lambda_1, \ldots, \lambda_k)$. The sum of the parts~$n$ is the \textit{weight} $\wt(\lambda)$ of $\lambda$. 

Let $\mathcal{P}$ denote the set of all partitions, including the empty partition of $n=0$. The enumeration of partitions by weight is equivalent to the classical question of counting the number of partitions of each non-negative integer. The resulting generating function 
\[Z_0(q) = \sum_{\lambda \in \mathcal{P}} q^{\wt(\lambda)} = 1 + 2q + 3q^2 + 5q^4 + 7 q^5+\ldots\]
has the well-known infinite product form \begin{equation}\label{eq:partition-generating-fct} Z_0(q)  = \prod_{k= 1}^{\infty}\frac{1}{1-q^{k}}.\end{equation}

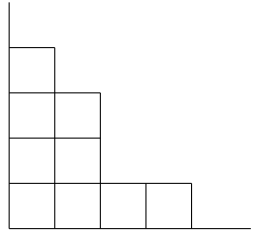
\begin{figure}\begin{center}
    \begin{tikzpicture}[scale=0.6, font=\footnotesize, fill=black!20]
	\draw (0, 0) -- (5.3,0);
	\draw (0,1) -- (4,1);
	\draw (0,2) -- (2,2);
	\draw (0,3) -- (2,3);
	\draw (0,4) -- (1,4);
	\draw (0,0) -- (0,5);
	\draw (1,0) -- (1,4);
	\draw (2,0) -- (2,3);
	\draw (3,0) -- (3,1);
	\draw (4,0) -- (4,1);
    \end{tikzpicture}
\caption{Young diagram of the partition $\lambda=(4,2,2,1)$ of $n=9$}
\label{fig:YD}
\end{center}\end{figure}

One can represent a partition $\lambda=(\lambda_1, \ldots, \lambda_k)$ in a planar arrangement, as a subset of the non-negative integer quadrant $\N\times \N$ as follows: one considers horizontal, left-adjusted bars of length $\lambda_1, \lambda_2, \ldots$, forming a left-and-bottom-adjusted subset of $\N\times\N$. We will refer to such subsets as Young diagrams, identifying partitions and Young diagrams below. A Young diagram consists of a finite union of {\em boxes} (sometimes also called {\em blocks}) $(a_i,b_i)\in \N\times\N$, with $a_i, b_i\geq 0$. See Figure~\ref{fig:YD} for an example. 

We next introduce a labelling scheme of our boxes $(a,b)\in \N\times \N$ which we call the {\it pattern of type~$A_r$}; this will be motivated by ideas from representation theory in the next section. Fix a positive integer $r$, and let $I=\Z/(r+1)$. We label each element $(a,b)\in\N\times \N$ by the label 
\[[(a-b)\!\!\!\!\mod (r+1)]\in I.\] This gives a periodic, diagonal labeling of $\N\times \N$, illustrated on Figure~\ref{fig:pattern-A}. 

\begin{figure}
    \begin{tikzpicture}[scale=0.6, font=\footnotesize, fill=black!20]
    \foreach \x in {1,2,4,5,6,7,8}{\draw (\x, 0) -- (\x,6.3);}
    \draw (0, 0) -- (0,7);
    
    \foreach \y in {1,2,4,5,6}{\draw (0,\y) -- (8.3,\y);}
    \draw (0,0) -- (9.3,0);
    
    \draw (0.5,0.5) node {0};
    \draw (1.5,0.5) node {1};
    \draw (4.5,0.5) node {$r$$-$$1$};
    \draw (5.5,0.5) node {$r$};
    \draw (0.5,1.5) node {$r$};
    \draw (1.5,1.5) node {0};
    \draw (4.5,1.5) node {$r$$-$$2$};
    \draw (5.5,1.5) node {$r$$-$$1$};
    \draw (6.5,0.5) node {0};
    \draw (7.5,0.5) node {1};
    \draw (6.5,1.5) node {$r$};
    \draw (7.5,1.5) node {0};
    
    \draw (0.5,4.5) node {1};
    \draw (1.5,4.5) node {2};
    \draw (0.5,5.5) node {0};
    \draw (1.5,5.5) node {1};
    \draw(0.5,3.1) node {\vdots};
    \draw(1.5,3.1) node {\vdots};
    \draw(3,0.5) node {\dots};
    \draw(3,1.5) node {\dots};
    \draw(8.75,3) node {\dots};
    \draw(3,6.75) node {\vdots};
    \end{tikzpicture}    
    \caption{The pattern of type $A_r$: periodic labelling of $\N\times \N$ with $(r+1)$ labels}
    \label{fig:pattern-A}
\end{figure}

\begin{figure}
    \begin{tikzpicture}[scale=0.6, font=\footnotesize, fill=black!20]

    \filldraw[draw=lightgray,fill=lightgray] (0,4) -- (0,2) -- (2,2) -- (2,3) -- (1,3) -- (1,4);
    
	\draw (0, 0) -- (5.3,0);
	\draw (0,1) -- (4,1);
	\draw (0,2) -- (2,2);
	\draw (0,3) -- (2,3);
	\draw (0,4) -- (1,4);
	\draw (0,0) -- (0,5);
	\draw (1,0) -- (1,4);
	\draw (2,0) -- (2,3);
	\draw (3,0) -- (3,1);
	\draw (4,0) -- (4,1);
	\draw (0.5,0.5) node {0};
	\draw (1.5,0.5) node {1};
	\draw (2.5,0.5) node {2};
	\draw (0.5,1.5) node {2};
	\draw (1.5,1.5) node {0};
	\draw (0.5,2.5) node {1};
	\draw (1.5,2.5) node {2};
	\draw (0.5,3.5) node {0};
	\draw (3.5,0.5) node {0};
    \end{tikzpicture}
    
    \caption{Young diagram of the partition $\lambda=(4,2,2,1)$ labelled by $I=\Z/3$, of multiweight (4,2,3). The gray boxes form the only removable border strip of length~$3$}
    \label{fig:first-young-diagram}
\end{figure}

Given a partition $\lambda\in\CP$, represent it by its Young diagram, a subset of $\N\times \N$. The \textit{multi-weight} of $\lambda$ is defined to be the vector
$(\wt_0(\lambda), \ldots, \wt_r(\lambda))$, where $\wt_i(\lambda)$ counts the number of boxes in~$\lambda$ of label $i\in I$. Consider the generating function
\begin{equation} \label{formulaZr} Z_r(q_0, \ldots, q_r) = \sum_{\lambda \in \mathcal{P}} \prod_{i\in I} q_i^{\wt_i(\lambda)}
\end{equation}
counting Young diagrams in $\Part$ by their full multiweight. 

To find a closed form for this generating function, one can make use of the so-called \textit{Littlewood decomposition}. The starting point of this correspondence is the following procedure. Given a partition $\lambda$, a \textit{border strip} of length $(r+1)$ is an edge-connected set of boxes in (the Young diagram of) $\lambda$, each furthest away from the axes along its diagonal (i.e.~contained in the top-and-right boundary of $\lambda$). A border strip is considered \textit{removable}, if removing its boxes from $\lambda$ produces another Young diagram $\lambda'$. $\lambda$ is called an \textit{$(r+1)$-core} partition, if the corresponding diagram has no removable border strips. Note that by definition, a border strip of length $(r+1)$ contains exactly one box of each label $i\in I$.

The Littlewood decomposition describes how any Young diagram can be obtained by adding border strips to a core diagram. Denote by 
$\CorePart{r+1}$ the set of $(r+1)$-core partitions. 

\begin{theorem}[Littlewood decomposition~\cite{Littlewood}]\label{thm:littlewood}
    There is a bijection $$\Part \leftrightarrow \Part^{r+1} \times \CorePart{r+1} \; $$ defined by the following two properties:
    \begin{enumerate} 
    \item the projection $\Part \to \CorePart{r}$ maps any Young diagram to the unique core diagram obtained through removing boundary strips of length $(r+1)$ in any order; 
    \item a partition $\mu_i$ in the $i$-th copy of $\Part$ on the right hand side corresponds to adding $\wt(\mu_i)$ boundary strips of length $r+1$ with first label $i$ to an $(r+1)$-core diagram.
    \end{enumerate}
\end{theorem}
The fact that such a bijection exists is by no means obvious; even the fact that removing border strips from a single starting diagram in some arbitrary order should result in a unique final $(r+1)$-core is not obvious. An attractive combinatorial proof can be given using the language of Maya diagrams, see e.g.~\cite{JamesKerber, gyenge2018euler}. The partitions $\mu_i$ appearing in the first factor under the correspondence are sometimes refereed to as the {\it $r$-quotients} of~$\lambda$. 

\begin{example} \label{ex:checker} For $r=1$, the labelling of $\N\times \N$ is given by the checkerboard (alternating) labelling with two labels. A border strip of length $2$ is simply a domino shape, with one box of each label. The $2$-core diagrams are the staircase diagrams corresponding to the partitions $\lambda=(k,k-1,\ldots, 1)$. 
\end{example}

Consider the generating function
$$Z_{r, \text{core}}(q_0, \ldots, q_r)=  \sum_{\lambda \in \CorePart{r+1}} \prod_{i\in I} q_i^{\wt_i(\lambda)}$$
of $(r+1)$-core partitions. 
\begin{theorem}[James--Kerber~\cite{JamesKerber}] \ 
\label{thm:gen-fct-core}
\begin{enumerate}
\item There exists a combinatorial bijection
\[ \CorePart{r+1} \longleftrightarrow \Z^r.\] 
\item     The generating function of  $(r+1)$-core partitions is given by 
    $$Z_{r, {\rm core}}(q_0, \ldots, q_r)= \sum_{w \in \Z^r}  q^{\frac{1}{2} w^{\top} C_{A_r} w}\prod_{i=1}^r q_i^{w_i} \; ,$$ where $q = \prod_{i \in I} q_i$, and $$C_{A_r} = \begin{pmatrix} 2 & -1 & 0 & \\ -1 & 2 & \ddots & 0 \\ 0 & \ddots & \ddots & -1 \\ & 0 & -1 & 2 \end{pmatrix}$$ is the $ r\times r$ Cartan matrix of {\em finite} type $A_r$.
    \end{enumerate}
\end{theorem}

Notice that even though this entire section has been of a combinatorial nature, the right Dynkin diagram has once again entered the picture through the Cartan matrix in this statement. 

We can combine this result with Theorem \ref{thm:littlewood} to obtain the full generating function $Z_r$: adding a border strip of length $(r+1)$ to any Young diagram adds exactly one box of each label. In the generating function, this corresponds to multiplying an entry by $q = \prod_{i \in I} q_i$. Using also~\eqref{eq:partition-generating-fct}, we obtain
\begin{corollary}\label{cor_Zr} The generating function \eqref{formulaZr} of all partitions labelled in the pattern of type $A_r$ can be expressed as
    \begin{equation}\label{eq:Zr} Z_r(q_0, \ldots, q_r) = (Z_0(q))^{r+1}\cdot Z_{r, {\rm core}}(q_0, \ldots, q_r) = \frac{\displaystyle\sum_{w \in \Z^r}  q^{\frac{1}{2} w^{\top} C_{A_r} w}\prod_{i=1}^r q_i^{w_i}}{\displaystyle\prod_{k = 0}^{\infty} (1-q^k)^{r+1}} \; ,\end{equation} where $q = \prod_{i \in I} q_i$.
\end{corollary}
This formula was re-proved in a completely elementary, combinatorial way in~\cite{GyengeEnumeration}. 
\begin{example}\label{ex:checker2} Continuing with the example $r=1$ from Example~\ref{ex:checker}, the generating function $Z_1(q_0, q_1)$ counts partitions labelled (coloured) as the checkerboard. We get the generating function 
\begin{equation} \label{eq:checker} Z_1(q_0, q_1) = Z_0(q_0q_1)^2 \cdot \sum_{m=-\infty}^{\infty} q_0^{m^2} q_1^{m^2+m} = \prod_{k=1}^\infty \frac{(1 + q_0^{2k-1}q_1^{2k})(1 + q_0^{2k-1}q_1^{2k-2})}{(1-q_0^kq_1^k)(1-q_0^{2k-1}q_1^{2k-1})},
\end{equation}
where the second equality uses a form of the Jacobi triple product identity. 
\end{example}

We close this section with a generalisation of Example~\ref{ex:checker}-\ref{ex:checker2}. Consider the labelling scheme where we label a box $(a,b)\in \N\times \N$
by the pair $(a\mbox{ mod } 2 \mid b \mbox{ mod } 2)$, another periodic labelling of $\N\times \N$ illustrated on Figure~\ref{fig:Z2Z2}. This labelling appears in~\cite{Boulet}. Given a partition $\lambda$, we let $\wt_{00}(\lambda), \ldots, \wt_{11}(\lambda)$ count boxes in (the Young diagram of) $\lambda$ with the appropriate labels.

\begin{theorem}[Boulet~\cite{Boulet}] \label{thm:fourparameter}
The generating function 
$$Z_{1,1}(q_{00}, q_{01},q_{10}, q_{11})=  \sum_{\lambda \in \CP} \prod_{i,j=0}^1 q_{ij}^{\wt_{ij}(\lambda)}$$
admits the infinite product form
\begin{equation}\label{eq:fourparameter}
Z_{1,1}(q_{00}, q_{01},q_{10}, q_{11})= \prod_{k=1}^\infty \frac{(1+q_{00}^kq_{10}^{k-1}q_{01}^{k-1}q_{11}^{k-1})(1+q_{00}^kq_{10}^kq_{01}^kq_{11}^{k-1})}{(1-q_{00}^kq_{10}^{k}q_{01}^{k}q_{11}^{k})(1-q_{00}^kq_{10}^{k-1}q_{01}^{k}q_{11}^{k-1})(1-q_{00}^kq_{10}^{k}q_{01}^{k-1}q_{11}^{k-1})}.
\end{equation}
\end{theorem}
Upon specialisation $q_{00}=q_{11}=q_0$ and $q_{01}=q_{10}=q_1$, this result recovers the second formula in~\eqref{eq:checker}.

Note finally that the pattern of type $A_r$ considered before and this last labelling admit a common generalisation, where we label a box $(a,b)\in \N\times \N$ by the pair $(a\mbox{ mod } (r+1) \mid b \mbox{ mod } (r+1))$. The corresponding generating function $Z_{r,r}(q_{00}, \ldots, q_{rr})$, in $(r+1)^2$ 
variables, does not appear to have been considered in the literature before. The obvious specialisation of any formula for this generating function 
would have to return the formula from Corollary~\ref{cor_Zr}. 

\begin{figure}
    \begin{tikzpicture}[scale=0.6, font=\footnotesize, fill=black!20]
    \foreach \x in {1,2,3}{\draw (\x, 0) -- (\x,3.3);}
    \draw (0, 0) -- (0,4);
    
    \foreach \y in {1,2,3}{\draw (0,\y) -- (3.3,\y);}
    \draw (0,0) -- (4,0);
    
    \draw (0.5,0.46) node {$0|0$};
    \draw (1.5,0.46) node {$1|0$};
    \draw (2.5,0.46) node {$0|0$};
    \draw (0.5,1.46) node {$0|1$};
    \draw (1.5,1.46) node {$1|1$};
    \draw (2.5,1.46) node {$0|1$};
    \draw (0.5,2.46) node {$0|0$};
    \draw (1.5,2.46) node {$1|0$};
    \draw (2.5,2.46) node {$0|0$};
    
    \draw(1.5,3.7) node {\vdots};
    \draw(3.7,1.5) node {\dots};
    
    \end{tikzpicture}    
    \caption{Periodic labelling of $\N\times \N$ with a pair of labels mod~$2$}
    \label{fig:Z2Z2}
\end{figure}

\subsection{Equivariant Hilbert schemes and labelled partitions}\label{subsec:HilbandP} In this section, we relate the combinatorics of the previous section to Hilbert schemes. The starting point is well known. Recall that the Hilbert scheme $\Hilb(\Aff^2)$ of the affine plane $Z =\Aff^2$ parametrises finite-dimensional quotients $R\twoheadrightarrow R/\ideal$ of the coordinate ring $R=\C[\Aff^2]\cong \C[x,y]$, corresponding to finite-colength ideals $\ideal\lhd R$. Among these ideals, there is a set of distinguished ones, defined by the condition that~$\ideal$ is generated by a finite collection of monomials $x^ay^b$. These ideals are exactly the fixed points on $\Hilb(\Aff^2)$ of the natural action of the algebraic torus 
$T=(\C^*)^2$ that rescales the coordinates $(x,y)$. By associating to a monomial ideal $\ideal\lhd R$ the set of boxes $(a,b)\subset\N\times\N$ with the condition that $x^ay^b\not\in \ideal$, we obtain a partition $\lambda$, the weight of which equals the codimension of $\ideal$. See Figure~\ref{fig:idealpartition} for an example. 

\begin{figure}\begin{center}
    \begin{tikzpicture}[scale=0.6, font=\footnotesize, fill=black!20]
	\draw (0, 0) -- (5.3,0);
	\draw (0,1) -- (4,1);
	\draw (0,2) -- (2,2);
	\draw (0,3) -- (2,3);
	\draw (0,4) -- (1,4);
	\draw (0,0) -- (0,5);
	\draw (1,0) -- (1,4);
	\draw (2,0) -- (2,3);
	\draw (3,0) -- (3,1);
	\draw (4,0) -- (4,1);
    \draw (0.5,0.5) node {$\phantom{1}1^{\phantom{1}}$};
    \draw (1.5,0.5) node {$\phantom{1}x^{\phantom{1}}$};
    \draw (2.5,0.5) node {$x^2$};
    \draw (3.5,0.5) node {$x^3$};
    \draw (4.5,0.5) node {{\color{red}$x^4$}};

    \draw (0.5,1.5) node {$\phantom{1}y^{\phantom{1}}$};
    \draw (0.5,2.5) node {$y^2$};
    \draw (0.5,3.5) node {$y^3$};
     \draw (0.5,4.5) node {{\color{red}$y^4$}};
    \draw (1.5,1.5) node {$\phantom{1}xy^{\phantom{1}}$};
    \draw (1.5,2.5) node {$xy^2$};
     \draw (2.5,1.5) node {{\color{red}$x^2y$}};
     \draw (1.5,3.5) node {{\color{red}$xy^3$}};
    \end{tikzpicture}
\caption{The partition $\lambda=(4,2,2,1)$ corresponds to the monomial ideal
$\ideal=\langle x^4, x^2y,xy^3,y^4\rangle$}
\label{fig:idealpartition}
\end{center}\end{figure}

We get the following result, a version of the famous G\"ottsche formula. 
\newcommand{\tp}{{\mathrm{top}}}
\begin{theorem} \label{thm_r0}\begin{enumerate} 
\item There is a one-to-one correspondence between $T$-fixed points on the Hilbert scheme $\Hilb(\Aff^2)$ and the set of partitions $\CP$. 
\item For each $n\geq 0$, the Hilbert scheme $\Hilb^n(\Aff^2)$ is (nonempty and) irreducible. The generating function of topological Euler characteristics satisfies
\begin{equation} \label{eq:Hilbert:char} \sum_{n\in \N} \chi_\tp(\Hilb^n(\Aff^2)) q^n = Z_0(q),\end{equation}
where $Z_0(q)$ is the generating function~\eqref{eq:partition-generating-fct} of partitions. 
\end{enumerate}\end{theorem}

We want to generalise this discussion to the situation of an abelian group acting on $\Aff^2$. Recall that the finite subgroup of type $A_r$ of $\SL(2,\C)$ is the cyclic group $$G_{A_r} = \left\langle \begin{pmatrix} \zeta & 0 \\ 0 & \zeta^{-1} \end{pmatrix}\right\rangle < \SL(2, \C) \; $$ of order $r+1$, where $\zeta$
% = \exp \left( \frac{2 \pi \sqrt{-1}}{r+1}\right)
is a primitive $(r+1)$-st root of unity; write $\sigma\in G_{A_r}$ for its given generator. The irreducible representations $\rho_i$ of $G_{A_r}$ are all one-dimensional and %in one-to-one correspondence with the set of $(r+1)$-st roots of unity: $V_i = \C$ as a vector space, and $\sigma(v) = \zeta^i v$ for $v \in V_i$. Under tensor product, the irreducible representations $\{\rho_i\}$ 
form the dual group $\Z/(r+1)$ under tensor product, which we will identify with the set~$I$ of labels from the previous section. 

The key fact we are going to rely on is that the action of~$G_{A_r}$ and the action of~$T$ on~$\Aff^2$ commute. Hence the torus $T=(\C^*)^2$ acts on the equivariant Hilbert scheme $\Hilb(\Aff^2)^{G_{A_r}}$. 

Consider the action of $G_{A_r}$ on $R = \C[x,y]$. Every monomial $x^a y^b$ is an eigenvector with respect to this action, decomposing $R$ into irreducible representations. The representation type of the subspace spanned by a monomial $x^a y^b$ is $\rho_i$, where $i=[(a-b)\!\!\mod (r+1)]\in I$. We thus recover the labelling of $\N\times\N$ considered in the previous section, displayed in Figure~\ref{fig:pattern-A}.

Further, monomial ideals $\ideal\lhd R$ of finite colength are automatically $G_{A_r}$-equivariant, and so are in the equivariant Hilbert scheme $\Hilb(\Aff^2)^{G_{A_r}}$. The following result generalises the non-equivariant Theorem~\ref{thm_r0}. 
\begin{theorem}[\cite{Gordon, Fujii-Minabe}] \label{thm:typeAmain}
Consider the equivariant Hilbert scheme \[\Hilb(\Aff^2)^{G_{A_r}}= \bigsqcup_{\rho\in\Rep(G_{A_r})} \Hilb^\rho(\Aff^2)\] attached to the cyclic group $G_{A_r}<\SL(2,\C)$. 
\begin{enumerate}\item For $\rho\in\Rep(G_{A_r})$, the space $\Hilb^\rho(\Aff^2)$, if non-empty, is irreducible. The set of components of 
$\Hilb(\Aff^2)^{G_{A_r}}$ is in bijection with the set $\CorePart{r}\times \N$ consisting of pairs $(\lambda, n)$ of an $(r+1)$-core partition $\lambda\in \CorePart{r}$ and a non-negative integer $n\in\N$. 
\item The $(\C^*)^2$-fixed points on $\Hilb(\Aff^2)^{G_{A_r}}$ are isolated, and are in one-to-one correspondence with $I$-labelled partitions $\lambda\in\CP$. For an $I$-labelled partition $\lambda$, the $G_{A_r}$-representation of the corresponding quotient $R/\ideal$ is isomorphic to $\oplus_{i\in I} \rho_i^{\wt_i(\lambda)}$. The point $[R/\ideal]\in \Hilb(\Aff^2)^{G_{A_r}}$ is contained in the component corresponding to the pair $(\mu, n)$ with $\mu$ being the $(r+1)$-core of $\lambda$, and $n$ the total weight of its $(r+1)$-quotients. 
\item We have an equality
\begin{equation}\label{eq:Zrchar_topEuler}\sum_{(n_i)\in \N^I} \chi_\tp\left(\Hilb^{\oplus_{i\in I} \rho_i^{n_i}}(\Aff^2)\right) \prod_{i\in I}q_i^{n_i} = Z_r(q_0, \ldots, q_r)
\end{equation}
between the generating function of topological Euler characteristics of components of $\Hilb(\Aff^2)^{G_{A_r}}$ and the generating function $Z_r$ of $I$-labelled partitions defined in~\ref{formulaZr}, written in a closed form in Corollary~\ref{cor_Zr}. 
\end{enumerate}
\end{theorem}

\subsection{Some abelian groups in $\GL(2,\C)$}\label{sec:GL2} Let us briefly comment on an extension of the ideas of the previous sections to abelian groups $G<\GL(2,\C)$ that are not contained in $\SL(2,\C)$. The first natural candidate is the subgroup 
$$G_{r+1,a} = \left\langle \begin{pmatrix} \zeta & 0 \\ 0 & \zeta^{a} \end{pmatrix}\right\rangle < \GL(2, \C) \; $$
where $\zeta$ is a primitive $(r+1)$st root of unity, and $1<a<r-2$. This group is abstractly still isomorphic to $C_{r+1}$. There is no problem with the definition of the equivariant Hilbert scheme $\Hilb(\Aff^2)^{G_{r+1,a}}$, the torus action, and the identification of the fixed points using the obvious generalisation of the labelling scheme. However, the resulting generating function appears much more difficult to study. A combinatorial approach was presented in~\cite{gyengePMH} and some periodicity results proven in~\cite{BZ}, but no closed form is known. 

Let us finally return briefly to the labelling scheme appearing before the statement of Theorem~\ref{thm:fourparameter}. It is easy to see that this labelling corresponds to the action of the group 
\[C_2\times C_2 \cong H=\left(\begin{array}{cc} \pm 1 & 0  \\0 & \pm 1\end{array}\right)<\GL(2,\C)\] 
on the affine plane $\Aff^2$, and thus on the Hilbert scheme $\Hilb(\Aff^2)$. It is interesting that even in this, non-$\SL$-case, an evaluation of the generating function of Euler characteristics of the equivariant Hilbert scheme $\Hilb(\Aff^2)^H$ in the closed form~\eqref{eq:fourparameter} is possible. The more general question asked at the end of~\ref{subsec:part} corresponds to the analogous action of the group $C_{r+1}\times C_{r+1}<\GL(2,\C)$; a complete solution to that problem would specialise to a solution of the problem raised in the previous paragraph. 

\subsection{Quot schemes and specialised generating functions}

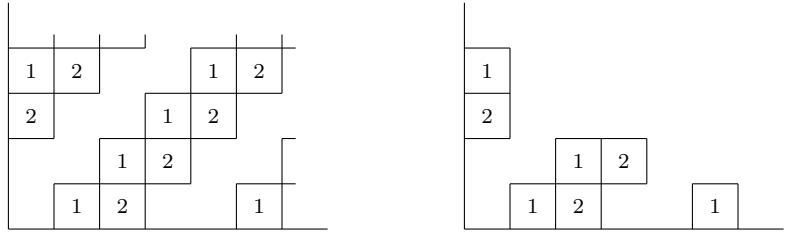
\begin{figure}
    \begin{tikzpicture}[scale=0.6, font=\footnotesize, fill=black!20]
	\draw (0,0) -- (7,0);
	\draw (0,0) -- (0,5);
    \draw (1,0) -- (1,1);
	\draw (2,0) -- (2,2);
	\draw (3,0) -- (3,3);
	\draw (4,1) -- (4,4);
	\draw (5,2) -- (5,4.3);
	\draw (6,3) -- (6,4.3);
    \draw (1,1) -- (4,1);
	\draw (2,2) -- (5,2);
    \draw (3,3) -- (6,3);
	\draw (4,4) -- (6.3,4);

    \draw (0,2) -- (1,2);
	\draw (0,3) -- (2,3);
	\draw (0,4) -- (3,4);
	\draw (1,2) -- (1,4.3);
	\draw (2,3) -- (2,4.3);
	\draw (3,4) -- (3,4.3);

    \draw (5,0) -- (5,1);
	\draw (6,0) -- (6,2);
	\draw (5,1) -- (6.3,1);
	\draw (6,2) -- (6.3,2);
 
	\draw (1.5,0.5) node {1};
	\draw (2.5,0.5) node {2};
	\draw (2.5,1.5) node {1};
	\draw (3.5,1.5) node {2};
	\draw (3.5,2.5) node {1};
	\draw (4.5,2.5) node {2};
	\draw (4.5,3.5) node {1};
	\draw (5.5,3.5) node {2};

	\draw (0.5,2.5) node {2};
	\draw (0.5,3.5) node {1};
	\draw (1.5,3.5) node {2};

 	\draw (5.5,0.5) node {1};

%---------------------------------
    \draw (10,0) -- (17,0);
	\draw (10,0) -- (10,5);
    \draw (11,0) -- (11,1);
	\draw (12,0) -- (12,2);
	\draw (13,0) -- (13,2);
	\draw (14,1) -- (14,2);
    \draw (11,1) -- (14,1);
	\draw (12,2) -- (14,2);

    \draw (10,2) -- (11,2);
	\draw (10,3) -- (11,3);
	\draw (10,4) -- (11,4);
	\draw (11,2) -- (11,4);

    \draw (15,0) -- (15,1);
	\draw (16,0) -- (16,1);
	\draw (15,1) -- (16,1);
 
	\draw (11.5,0.5) node {1};
	\draw (12.5,0.5) node {2};
	\draw (12.5,1.5) node {1};
	\draw (13.5,1.5) node {2};

	\draw (10.5,2.5) node {2};
	\draw (10.5,3.5) node {1};

 	\draw (15.5,0.5) node {1};
    \end{tikzpicture}
    \caption{The pattern of type $A_{3,\{1,2\}}$ (left) and an example of a Young diagram in this pattern (right)}
    \label{fig:young-diagram-A-r-J}
\end{figure}

For a nonempty subset $J\subset I$, recall the Quot scheme $\Quot_{G_{A_r},J}(\Aff^2)$ from Definition-Theorem~\ref{defthm:quotscheme}.
The modules $e_J(R \rtimes G_{A_r})e_0$ appearing in the definition of this Quot scheme can be visualised in a manner similar to the case $J=I$ from the previous section. Consider the subset $(\N\times\N)_J$ of $\N\times\N$ consisting of all boxes with labels in~$J$. We call the resulting arrangements of labelled boxes \textit{the pattern of type $A_{r,J}$}. Write $\Part[r,J]$ for the set of finite subsets of $(\N\times\N)_J$ that can be obtained by intersecting the pattern of type $A_{r,J}$ with a Young diagram. Note that for $\nu\in\Part[r,J]$, the Young diagram $\lambda$ realising it is usually not unique. 
Examples of a pattern of type $A_{r,J}$ and an element of $\Part[r,J]$ are shown in Figure~\ref{fig:young-diagram-A-r-J}.

%Now while the Hilbert and Quot schemes associated to the singularity contain many more points than just the monomial quotient modules representable by Young diagrams, we can still extract information from these isolated points. The key tool here is that of a \textit{torus action}. The two-dimensional torus $(\C^*)^2$ acts on the affine plane $V$ by rescaling along the coordinate axes, and this action and commutes with the action of $G$ (this is only true in type $A$). This action furthermore gives rise to an action on the Hilbert and Quot schemes:

%The notion of monomial submodules and the corresponding quotients makes sense; these are in bijection with the following set: the Young diagram representing a monomial quotient of $e_J(R \rtimes G)e_0$

For each $\nu\in \Part[r,J]$ and $i\in J$, the notion of the $i$-weight $\wt_i(\nu)$ is well defined, and we can consider the generating function 
$$Z_{r,J}(q_i\colon i\in J) = \sum_{\nu \in \Part[r,J]} \prod_{i\in J}q_i^{\wt_{i}(\nu)} \; .$$ 

\begin{theorem}[\cite{BGS}] \label{thm:fixed-points-young-diagrams} For nonempty subsets $J\subset I$, consider the Quot schemes
\[\Quot_{G_{A_r},J}(\Aff^2) = \bigsqcup_{v_J \in \N^J} \Quot_{G_{A_r},J}^{v_J}(\Aff^2) \; .\]
\begin{enumerate} 
\item The torus $T=(\C^*)^2$ acts on each Quot scheme $\Quot_{G_{A_r},J}(\Aff^2)$ in such as way that the degeneration maps~\eqref{eq:quot-degeneration} are $(\C^*)^2$-equivariant. 
\item For a fixed non-empty subset $J\subset I$, the fixed points of the $T$-action on $\Quot_{G_{A_r},J}(\Aff^2)$ are the points parametrising monomial quotients, which are in bijection with the set $\Part[r,J]$ of Young diagrams in the pattern of type $A_{r,J}$.
\item We have an equality of generating functions
$$Z_{r,J}(q_j\colon j\in J) = \sum_{v_J \in \N^J} \chi_\tp(\Quot_{G_{A_r}, J}^{v_J}(\Aff^2))\prod_{j \in J} q_j^{v_j} \; .$$ 
\end{enumerate}
\end{theorem}

We are not aware of a straightforward combinatorial characterisation of the components of a Quot scheme similar 
to Theorem~\ref{thm:typeAmain}(1). 

%Torus actions are a powerful tool in understanding the geometry of a space. For example, when a torus acts on a variety with finitely many isolated fixed points, the Euler characteristic of the variety equals the number of fixed points. {\balazs Cite...} Therefore, and by Theorem \ref{thm:fixed-points-young-diagrams}, the question of calculating Euler characteristics of the Quot schemes is equivalent to enumerating Young diagrams.

% We then know that $\chi(\Quot_{G, J}^{v_J}(\Aff^2))$ equals the number of torus fixed points in $\Quot_{G, J}^{v_J}(\Aff^2)$, which are in bijective correspondence with all the Young diagrams in the pattern of type $A_{r,J}$ containing $v_j$ boxes labeled $j$ for each $j \in J$. 

%\subsubsection{Young diagrams in the pattern of type $A_{r,J}$ and the substitution formula}

To end this section, we explain a formula for the generating function $Z_{r,J}$, using a curious root-of-unity calculus from~\cite{gyenge2018euler, BGS}. Recall that $Z_{r,J}$ is a series in the variables $\{q_j: j \in J\}$. We will show that it can be obtained from the full generating function $Z_r$ from~\eqref{thm:typeAmain} through substituting certain roots of unity. 

In order to do this, we will make use of the degeneration maps (\ref{eq:quot-degeneration}) in the special case where the left-hand side is an equivariant Hilbert scheme. Recall from Theorem \ref{thm:fixed-points-young-diagrams}(1) that these maps are equivariant with respect to the torus action. As as a consequence, they map fixed points to fixed points. In terms of Young diagrams, this amounts to a map of sets $$\pi_{r,J} \colon \Part \to \Part[r,J] \; ,$$ given by intersecting any full Young diagram $\lambda\in\CP$ with the pattern of type $A_{r,J}$. The known counting function for $\Part$ can therefore be arranged along the fibers of this map, i.e.\ written as a sum over $\Part[r,J]$: $$Z_r(q_0, \ldots, q_r) = \sum_{\nu \in \Part[r,J]} \sum_{\substack{\lambda \in \Part \\ \pi_{r,J}(\lambda) = \nu}} \prod_{i\in I} q_i^{\wt_i(\lambda)} \; ,$$ We wish to relate this series to $$Z_{r,J}(q_i \colon i\in J) = \sum_{\nu \in \Part[r,J]} \prod_{i\in J} q_i^{\wt_{i}(\nu)} \; .$$ To do so, we will exhibit a substitution rule mapping the fibre sum $\sum_{\lambda \in \pi_{r,J}^{-1}(\nu)} \prod_{i\in I} q_i^{\wt_i(\lambda)}$ to $\prod_{i\in J} q_i^{\wt_i(\nu)}$ for each $\nu\in\Part[r,J]$.

\begin{figure}
    \begin{tikzpicture}[scale=0.3, font=\footnotesize, fill=black!20]
    \draw (0,0) -- (16.5,0);
    \draw (0,0) -- (0,14.5);

    \draw (1,12) -- (1,13);
    \draw (0,12) -- (1,12);
    \draw (0,13) -- (1,13);
    
    \draw (1,5) -- (1,8);
    \draw (2,6) -- (2,9);
    \draw (3,7) -- (3,9);
    \draw (4,8) -- (4,9);
    \draw (0,5) -- (1,5);
    \draw (0,6) -- (2,6);
    \draw (0,7) -- (3,7);
    \draw (1,8) -- (4,8);
    \draw (2,9) -- (4,9);
    
    \draw (1,0) -- (1,1);
    \draw (2,0) -- (2,2);
    \draw (3,0) -- (3,3);
    \draw (4,1) -- (4,4);
    \draw (5,2) -- (5,5);
    \draw (6,3) -- (6,6);
    \draw (7,4) -- (7,6);
    \draw (8,5) -- (8,6);
    \draw (1,1) -- (4,1);
    \draw (2,2) -- (5,2);
    \draw (3,3) -- (6,3);
    \draw (4,4) -- (7,4);
    \draw (5,5) -- (8,5);
    \draw (6,6) -- (8,6);

    \draw (8,0) -- (8,1);
    \draw (9,0) -- (9,2);
    \draw (10,0) -- (10,3);
    \draw (11,1) -- (11,4);
    \draw (12,2) -- (12,4);
    \draw (8,1) -- (11,1);
    \draw (9,2) -- (12,2);
    \draw (10,3) -- (12,3);
    \draw (11,4) -- (12,4);

    \filldraw[draw=black,fill=lightgray] (1,13) rectangle (3,9);
    \filldraw[draw=black,fill=lightgray] (4,9) rectangle (7,6);
    \filldraw[draw=black,fill=lightgray] (8,6) rectangle (12,4);
    \filldraw[draw=black,fill=lightgray] (12,3) rectangle (15,0);

    \draw (0.5,12.5) node {2};
    
    \draw (0.5,6.5) node {1};
    \draw (1.5,7.5) node {1};
    \draw (2.5,8.5) node {1};
    \draw (0.5,5.5) node {2};
    \draw (1.5,6.5) node {2};
    \draw (2.5,7.5) node {2};
    \draw (3.5,8.5) node {2};
    
    \draw (1.5,0.5) node {1};
    \draw (2.5,1.5) node {1};
    \draw (3.5,2.5) node {1};
    \draw (4.5,3.5) node {1};
    \draw (5.5,4.5) node {1};
    \draw (6.5,5.5) node {1};
    \draw (2.5,0.5) node {2};
    \draw (3.5,1.5) node {2};
    \draw (4.5,2.5) node {2};
    \draw (5.5,3.5) node {2};
    \draw (6.5,4.5) node {2};
    \draw (7.5,5.5) node {2};
    
    \draw (8.5,0.5) node {1};
    \draw (9.5,1.5) node {1};
    \draw (10.5,2.5) node {1};
    \draw (11.5,3.5) node {1};
    \draw (9.5,0.5) node {2};
    \draw (10.5,1.5) node {2};
    \draw (11.5,2.5) node {2};
    \end{tikzpicture}
    
    \caption{An element $\nu\in\Part[6,\{1,2\}]$ in the pattern of type $A_{6,\{1,2\}}$. The boundary of any Young diagram $\lambda$ in the fibre $\pi_{6,\{1,2\}}^{-1}(\nu)$ passes through the gray rectangles.}
    \label{fig:young-fibers}
\end{figure}
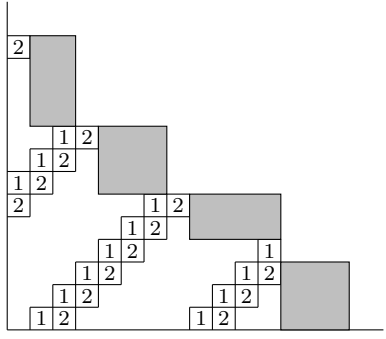

In order not to over-burden notation, we will explain the idea in an example. Let $r=6$ and $J=\{1,2\}$, a situation 
shown on Figure \ref{fig:young-fibers}. Fix a particular $\nu \in \Part[6,\{1,2\}]$ as on the Figure, and let us investigate the fibre $\pi_{r,J}^{-1}(\nu)$ in this case. The boundary of any Young diagram $\lambda$ mapping to $\nu$ under $\pi_{6,\{1,2\} }$ has to pass through the gray rectangles. In fact, the left and bottom edges of the rectangles, together with the boundary around the peaks of $\nu$, form a minimal Young diagram $\nu^{\vee} \in \pi_{r,J}^{-1}(\nu)$, and every other member of the fibre can be obtained by stacking Young diagrams contained in the rectangles on top of~$\nu^{\vee}$. Hence we can reduce counts along the fibre of $\pi_{r,J}$ to counting Young diagrams inside of rectangles. The following result is well-known.

\begin{proposition}\label{prop:part-box-q-binom}
    Consider the set $\PartRect{a,b}$ of Young diagrams inside a rectangle of height $a$ and width $b$. The generating function counting its elements by weight equals the $q$-binomial coefficient: $$\sum_{\lambda \in \PartRect{a,b}} q^{\wt(\lambda)} = \binom{a+b}{a}_q =\frac{(1-q^{a+b}) \cdots (1-q^{b+1})}{(1-q^a) \cdots (1-q)} \; .$$
\end{proposition}

Notice that independent of the particular choice of $\nu$, all the gray rectangles appearing in the fibers have the same number $a+b=6$, which is one more than the number of nodes in~$I$ lying between the nodes $2$ and $1$ in cyclic order. We make the following general observation.

\begin{lemma}\label{lem:subst-q-binom}
    Let $a,b \in \N$. If $\xi$ is a primitive $(a+b+1)$-th root of unity, then substituting $q = \xi$ into the $q$-binomial coefficient gives $$\binom{a+b}{a}_{\xi} = (-1)^k \xi^{-\frac{1}{2} a(a+1)}.$$
\end{lemma}

This $(a+b+1)$-th root of unity is the one we want to use for our substitution. In the example above, we substitute the seventh root of unity $\xi$ for the $(I \setminus J)$-indexed variables $q_3$, $q_4$, $q_5$, $q_6$ and $q_0$. We furthermore replace $q_1$ by $\xi q_1$ and $q_2$ by $\xi q_2$. For fixed $\nu \in \Part[r,J]$ this substitution replaces the fiber generating function $\sum_{\lambda \in \pi_{r,J}^{-1}(\nu)} \prod_{i\in I} q_i^{\wt_i(\lambda)}$ by the monomial $\prod_{j\in J} q_j^{\wt_j(\nu)}$ multiplied by a certain root of unity. A calculation using Lemma \ref{lem:subst-q-binom} now shows that this root of unity factor is in fact independent of~$\nu$. Hence the substitution maps the entire power series $Z_r$ to this root of unity times $Z_{r,J}$.

This calculation works for arbitrary $r$ and nonempty $J \subset I$. In general, the full sub-diagram supported on $I \setminus J$ will be a union of some finite type $A$ diagrams. For any vertex $i \in I \setminus J$ we let $r_i$ denote the number of vertices in its connected component of this sub-diagram; for $i \in J$ we set $r_i=0$. Then the the following substitution formula lets us calculate $Z_{r,J}$ from $Z_r$.

%For a connected component of this subdiagram containing $r$ vertices, we fix a primitive $(r+2)$-th root of unity $\xi$ and substitute all the $q_i$ corresponding to vertices in this connected component by $\xi$. We also multiply the two $J$-vertices bordering this connected component by $\xi$. 

%It turns out, that this substitution, together with some correctional roots of unity multiplied to the remaining variables, and one to the entire generating function, gives precisely what we want. For a node $i \in I$ let $k_i$ be the number of consecutive nodes not in $J$ adjacent to the right of $i$ and, and $l_i$ the same to the left of $i$. {\balazs Lukas please write this more clearly.} In other words, $k_i$ and $l_i$ are the $r_j$ for the components of $I \setminus J$ adjacent to $i$, and zero if there is none. Then the general generating function $Z_{r,J}$ can be determined as follows.

\begin{theorem}[\cite{BGS}]\label{thm:subst-J-young-walls}
    Given $J \subseteq I$, consider the following substitution of variables $q_I$:
    \begin{equation}\label{Subst}
        q_i \mapsto s_{r, J}(q_i) = \begin{cases} q_i \exp \left(\frac{2 \pi \sqrt{-1}}{2 + r_{i-1}} + \frac{2 \pi \sqrt{-1}}{2 + r_{i+1}} \right), & i \in J, \\ \exp \left(\frac{2 \pi \sqrt{-1}}{2 + r_i} \right), & i \notin J. \end{cases}
    \end{equation}
    Then there is a specific complex root of unity $c_{r,J}$ such that the substitution formula $$Z_{r,J}(q_i\colon i\in J) = c_{r,J} s_{r,J}(Z_r(q_0,\ldots, q_n)) $$
    %Let $m$ denote the minimum of $I$ in $\{0, \ldots, n\}$. Then $$c_{n,I} = (-1)^m \exp \left( 2 \pi \sqrt{-1} \frac{m(m+1)}{2(l_m + 2)}\right)$$
    holds. 
\end{theorem}
A special case of this formula was proved earlier in~\cite{DijkgraafS, gyenge2018euler}, when $J=\{0\}$. Geometrically, by Example~\ref{ex:quots}(2) the relevant Quot scheme in this case is the Hilbert scheme of points of the singular surface. 

\begin{remark} The substitution phenomenon of Theorem~\ref{thm:subst-J-young-walls} resembles in some respects the celebrated {\it cyclic sieving phenomenon}~\cite{CyclicSieving, Sagan}. As in our case, cyclic sieving consists of substituting carefully chosen roots of unity into generating functions, in the presence of cyclic group actions, to obtain counts of related combinatorial quantities. We have not been able to make this analogy precise. 
\end{remark}

\subsection{Representation theory around coloured partitions}\label{subset:repth}\ 

Define (fermionic) Fock space to be the infinite-dimensional vector space
\[ \BF = \bigoplus_{\lambda\in\CP} \C \vect{\lambda}
\]
generated by a basis indexed by the set $\CP$ of all partitions, that we will continue to view as Young diagrams. Call a box $s \in \lambda$ in a Young diagram $\lambda\in\CP$ {\it removable}, if $\lambda'=\lambda\setminus\{s\}\in\CP$; call a box $t \not\in \lambda$ {\it addable},  if  $\lambda'=\lambda\cup\{t\}\in\CP$. Note that a removable box is nothing but a removable border strip of length~$1$. Moreover, fixing as before the positive integer~$r$, for $\lambda\in\CP$ we denote by $h_c(\lambda)\in\{\pm1,0\}$ the difference between the number of addable and removable boxes of label $c\in \Z/(r+1)$. 
Finally, denote as before by $\wt_i(\lambda) \in \N$ the number of boxes labeled $i\in\Z/(r+1)$ in $\lambda$. 

The starting point of the relationship between partition combinatorics and the theory of infinite-dimensional Lie algebras is the following construction. 
Define four sets of operators on $\BF$, indexed by $c\in\Z/(r+1)$, as follows:
\[\begin{array}{c} e_c \vect{P}  = \displaystyle\sum_{\substack{{P'=P\setminus\{s\}} \\ {l(s)\equiv c \!\!\!\!\mod r+1}}}  \vect{P'}, \ \ \ \  
f_c  \vect{P}  =  \displaystyle\sum_{\substack{{P'=P\cup\{s\}} \\ {l(s)\equiv c \!\!\!\!\mod r+1}}} \vect{P'}, \\[2.5em] 
h_c \vect{P}  =  h_c(P)\vect{P}, \ \ \ \ 
d_c \vect{P} =  \wt_c(P)\vect{P}.\end{array}
\]
The operators $e_{c}, f_{c}, h_{c}, d_{c}$ for $c\in \Z/(r+1)$ act on $\BF$, since for every fixed $P\in\CP$, the number of terms in each sum is finite. 

\begin{example} The construction is already interesting in the case $r=0$, corresponding to the unlabelled case. It is a fundamental fact, easy to check directly, that with $r=0$, for every Young diagram $P\in\CP$ we have $h_0(P)=1$. In other words, there is always one more way to add a box to a Young diagram than to remove one. This means that the operator $h_0$ above reduces to the identify. The key commutation relation between the operators in this case is 
\[ [e_0,f_0] = h_0 = {\rm Id}_{\BF}. 
\]
The operator $d_1$, on the other hand, behaves as a grading operator
\[ [e_0,d_0] = e_0, \ \ \ \ \ [f_0,d_0] = -f_0,
\]
giving $\BF$ the structure of a graded space. In this case, the algebra $\langle e_0, f_0\rangle$ is nothing but the basic Heisenberg (or Weyl) algebra, defined by a ``raising'' operator $f_0$ and a ``lowering'' operator $e_0$, describing the relation between the position and momentum operators in basic quantum-mechanics. Using a further set of operators corresponding to adding and removing border strips of arbitrary length, one can in fact extend this representation to an action of the infinite-dimensional Heisenberg algebra ${\mathfrak{Heis}}$, see~\cite{kac1990infinite, Tingley}. The formula~\eqref{eq:partition-generating-fct} in turn becomes the graded character of the representation $\BF$ of ${\mathfrak{Heis}}$. 

The appearance of the infinite dimensional Heisenberg algebra ${\mathfrak{Heis}}$ in our story was explained by the following result. 
\begin{theorem}[Grojnowski~\cite{Groj}, Nakajima~\cite{Nak, NakajimaBook}] There is a graded isomorphism of ${\mathfrak{Heis}}$-representations
\[ \BF \cong H^*(\Hilb(\Aff^2), \C)
\]
between Fock space and the cohomology of the Hilbert scheme of points on $\Aff^2$, where the action of generators of  ${\mathfrak{Heis}}$ on the right hand side is given by geometrically defined operators.   
\end{theorem}
The equality of the graded characters of the two sides is~\eqref{eq:Hilbert:char}.
\end{example}

Turning to the case $r>0$, it is still not difficult to check by direct computation that the operators $e_{c}, f_{c}, h_{c}, d_{c}$ satisfy the basic relations
\[ [e_c,f_{c'}] = \delta_{c,c'} h_c, 
\]
as well as Serre-type commutation relations known from the theory of simple finite-dimensional Lie algebras, and grading relations. In this way, we get the algebras~\cite{kac1990infinite, Tingley}
\[ \langle e_{c}, f_{c}, h_{c} \colon c\in\Z/(r+1)\Z \rangle \cong \widehat{\mathfrak{sl}}'_{r+1},
\] 
the derived algebra of the affine Lie algebra attached to $\widehat A_r$, and the full affine Lie algebra 
\[ \langle d_{0}; e_{c}, f_{c}, h_{c} \colon c\in\Z/(r+1)\Z \rangle \cong \widehat{\mathfrak{sl}}_{r+1}.
\] 
Note that $\BF$ is reducible as an $\widehat{\mathfrak{sl}}_{r+1}$-module, though it becomes irreducible if one introduces a further set of operators forming a Heisenberg algebra, leading to a representation of the larger algebra $\widehat{\mathfrak{gl}}_{r+1}$. 

Similar to the case $r=0$, the formula~\eqref{eq:Zr} becomes a graded character of the representation $\BF$. One further has the following result.  
\begin{theorem}[Nakajima~\cite{Nakajima, NakajimaAffine}]\label{thm_reps_TypeA} We have a graded isomorphism of $\widehat{\mathfrak{gl}}_{r+1}$-modules 
\begin{equation}\label{eq:reps:r}\BF\cong \bigoplus_{\rho\in{\rm Rep}(G_{A_r})} H^*({\rm Hilb}^\rho(\Aff^2), \C)
\end{equation}
between Fock space and the cohomology of the equivariant Hilbert scheme, where the $\widehat{\mathfrak{gl}}_{r+1}$-action on the right hand side is constructed geometrically via correspondences.
\end{theorem}
The equality of graded characters is expressed by~\eqref{eq:Zrchar_topEuler}. 

\begin{remark}\label{rem:crystal} One can go further. The isomorphism~\eqref{eq:reps:r} can in fact be extended to an isomorphism of representations of {\em quantum} affine algebras~\cite{NakajimaAffine}. The main relevance to the present discussion is the fact that the basis of $\BF$ given by partitions is one model of the so-called crystal basis~\cite{kang2002book} of this representation. This observation will be crucial in the next section. 
\end{remark}

We finally mention the fact that, given a non-empty subset $J\in I$, we get a subalgebra inside $\widehat{\mathfrak{sl}}_{r+1}$, the algebra
\[ {\mathfrak s}_J=\langle e_{c}, f_{c}, h_{c} \colon c\in J \rangle < \widehat{\mathfrak{sl}}_{r+1}.\]
For example, for $J=I\setminus\{0\}$, we have 
\[ {\mathfrak s}_J\cong {\mathfrak{sl}}_{r+1}.\]
More generally, for a connected interval $J=\{a,\ldots, b\}$ of length $h$, such as the case $J=\{1,2\}$ of length $2$ above, we have
\[ {\mathfrak s}_J\cong {\mathfrak{sl}}_{h+1},\]
in particular for $J=\{0\}$, we get 
\[ {\mathfrak s}_J\cong {\mathfrak{sl}}_{2}.\]
In~\cite[Prop.4.10]{BGS}, we used a vector space spanned by partitions in rectangles, discussed above, to construct all irreducible fundamental representations of these semisimple Lie algebras.

\section{The non-abelian cases}\label{sec:nonabelian}

\subsection{Type D: the challenge} 
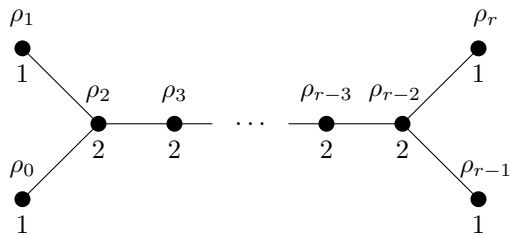
\begin{figure}[h]
 \begin{tikzpicture}
%    \foreach \i in {1,...,5}
%        \fill (\i, 0) circle(3pt);
    \fill (0,0) circle(3pt);
    \fill (0,2) circle(3pt);
    \fill (1,1) circle(3pt);
    \fill (2,1) circle(3pt);
    \fill (4,1) circle(3pt);
    \fill (5,1) circle(3pt);
    \fill (6,0) circle(3pt);
    \fill (6,2) circle(3pt);
       
    \draw (1,1) -- (2.5,1);
    \draw (3.5,1) -- (5,1);
    \draw (0,0) -- (1,1);
    \draw (0,2) -- (1,1);
    \draw (6,0) -- (5,1);
    \draw (6,2) -- (5,1);
    
    \node[below] at (0,-0.1) {$1$};
    \node[below] at (0,1.9) {$1$};
    \node[below] at (1, .9) {$2$};
    \node[below] at (2, .9) {$2$};
    \node[below] at (4, .9) {$2$};
    \node[below] at (5, .9) {$2$};
    \node[below] at (6, 1.9) {$1$};
    \node[below] at (6, -0.1) {$1$};
    
    \node at (0, .4) {$\rho_0$};
    \node at (0, 2.4) {$\rho_1$};
     \node at (1, 1.4) {$\rho_2$};
    \node at (2, 1.4) {$\rho_3$};
   \node at (4, 1.4) {$\rho_{r-3}$};
    \node at (4.9, 1.4) {$\rho_{r-2}$};
  \node at (3, 1) {$\ldots$};
    \node at (6.1, .4) {$\rho_{r-1}$};
    \node at (6.1, 2.4) {$\rho_r$};
    
\end{tikzpicture}
\caption{Labelling representations with their dimensions in the diagram of affine type $D_r$}
\label{fig:Dn_multiplicities}
\end{figure}

Let us replace the cyclic group $C_{r+1}\cong G_{A_r}< \SL(2,\C)$ of type~$A_r$ with the binary dihedral group 
$\tilde D_{2r-4}\cong G_{D_r} < \SL(2,\C)$ of type $D_r$, the case Theorem~\ref{thm:classification}(2). This group has $4$ one-dimensional representations $\rho_0$ (trivial), $\rho_1$ (sign), $\rho_{r-1}$ and $\rho_r$ corresponding to nodes at the ends of the affine Dynkin diagram of type $D_r$, and $(r-3)$ two-dimensional irreducible representations $\rho_2, \ldots, \rho_{r-2}$ corresponding to the intermediate nodes, with~$\rho=\rho_2$ being the representation corresponding to our standard two-dimensional vector space~$V\cong V_2$; see Figure~\ref{fig:Dn_multiplicities}. 

Attempting to follow the line of reasoning from the previous section, the immediate challenge one faces is that the action of the torus $T=(\C^*)^2$ on $\Aff^2$ no longer commutes with the action of $G_{D_r}$, only that of the constant diagonal torus $T_1\cong\C^*$ does. This smaller torus however no longer gives isolated fixed points on the equivariant Hilbert scheme, so 
the fundamental results Theorem~\ref{thm_r0}(1)-Theorem~\ref{thm:typeAmain}(2), giving an immediate connection to combinatorial ideas, have no direct analogue. 

The key to make progress comes from Remark~\ref{rem:crystal}, where we gave a representation-theoretic explanation for the appearance of labelled partitions. The important fact is that a version of the representation-theoretic Theorem~\ref{thm_reps_TypeA} holds in all types (see~\cite[Section 5]{GyNSzAnnouncement} for a more precise discussion). One then needs to find an analogous crystal basis construction in type $D_r$, and see whether the combinatorics can be related to Hilbert schemes. 

\subsection{Type D: Young walls}

We describe here the type $D$ analogue of diagonally colored partitions, following the crystal basis literature~\cite{kang2004crystal,kwon2006affine}. In this section, we formulate the construction as found in~[ibid.]; in the next section, we will explain how the construction fits into our context. 

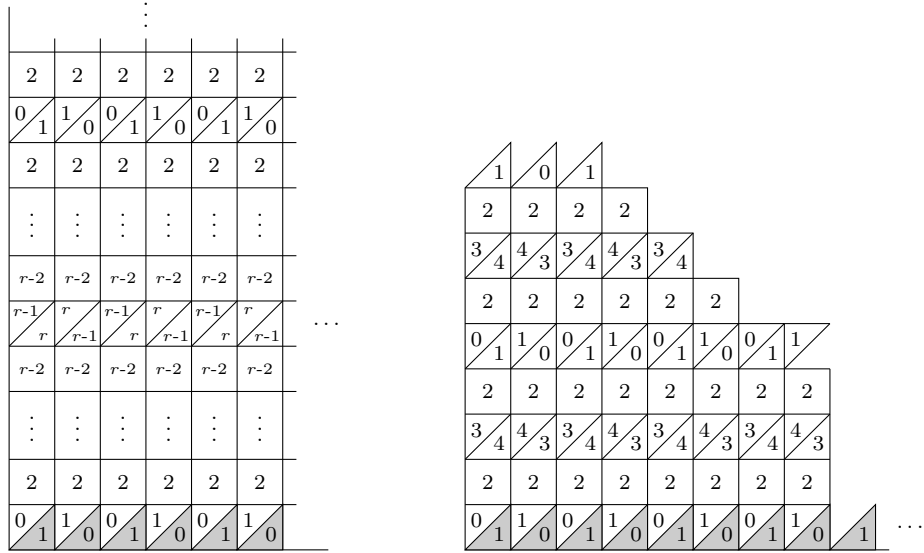
\begin{figure}[h]
    \begin{tikzpicture}[scale=0.6, font=\footnotesize, fill=black!20]
    \foreach \x in {1,2,3,4,5,6} {\draw (\x, 0) -- (\x,11.3);}
    \draw (0, 0) -- (0,12);
    
    \foreach \y in {1,2,3.5,4.5,5.5,6.5,8,9,10,11}{\draw (0,\y) -- (6.3,\y);}
    \draw (0,0) -- (7,0);
    
    \foreach \x in {0,1,2,3,4,5}{
    	%\draw (\x,0) -- (\x+1,1);
    	\draw (\x,4.5) -- (\x+1,5.5);
    	\draw (\x,9) -- (\x+1,10);
    	\draw (\x+0.5,1.5) node {2};
   
%     \draw (\x+0.5,4) node {$r$$-$$2$};
 %   	\draw (\x+0.5,6) node {$r$$-$$2$};
     \draw (\x+0.5,4) node {{\tiny $r$-$2$}};
    	\draw (\x+0.5,6) node {{\tiny $r$-$2$}};
  
    	\draw (\x+0.5,8.5) node {2};
    	\draw (\x+0.5,10.5) node {2};
    	\draw(\x+0.5,2.85) node {\vdots};
    	\draw(\x+0.5,7.35) node {\vdots};
    	\filldraw (\x,0) -- (\x+1,1) -- (\x+1,0) -- cycle;
     }
     \foreach \x in {0,2,4}{
        \draw (\x+0.25,0.65) node {0};
        \draw (\x+0.75,0.35) node {1};
%        \draw (\x+0.49,5.28) node  {$r$$-$$1$};
        \draw (\x+0.37,5.28) node  {{\tiny $r$-$1$}};
 
        	%\draw (\x+0.49,5.28) node %[fill=white] {$n$$-$$1$};
        \draw (\x+0.75,4.72) node {{\tiny $r$}};
        \draw (\x+0.25,9.65) node {0};
        \draw (\x+0.75,9.35) node {1};
    }
    \foreach \x in {1,3,5}{
        \draw (\x+0.25,0.65) node {1};
        \draw (\x+0.75,0.35) node {0};
        \draw (\x+0.25,5.28) node {{\tiny $r$}};
               \draw (\x+0.66,4.72) node {{\tiny $r$-$1$}};
 %       \draw (\x+0.52,4.72) node {$r$$-$$1$};
        \draw (\x+0.25,9.65) node {1};
        \draw (\x+0.75,9.35) node {0};
    }
    \draw (7,5) node {\dots};
    \draw (3,12) node {\vdots};

%----------------------------

    \draw (10, 0) -- (10,8) -- (11,9)--(11,8);
    \draw (11, 0) -- (11,8) -- (12,9)--(12,8);
    \draw (12, 0) -- (12,8) -- (13,9)--(13,8);
    \draw (13, 0) -- (13,8);
    \draw (14, 0) -- (14,8) -- (10,8);
    \draw (15, 0) -- (15,7) -- (10,7);
    \draw (16, 0) -- (16,6) -- (10,6);
    \draw (17, 0) -- (17,5) -- (10,5);
    \draw (17, 5) -- (18,5);
    \draw (18, 0) -- (18,4);
    \draw (10, 4) -- (18,4);
    \draw (10, 3) -- (18,3);
    \draw (10, 2) -- (18,2);
    \draw (10, 1) -- (18,1);
    \draw (19, 0) -- (19.3,0);
    %\draw (2, 0) -- (2,10) -- (0,10);
    %\draw (3, 0) -- (3,9) -- (0,9);
    %\draw (0, 0) -- (0,12);

    \draw (10,0) -- (19,0);
    \foreach \x in {0,1,2,3,4,5,6,7}{
	   %\draw (\x,0) -- (\x+1,1);
	   \draw (\x+10,2) -- (\x+11,3);
	   \draw (\x+10,4) -- (\x+11,5);
	   \draw (\x+10.5,1.5) node {2};
	   \draw (\x+10.5,3.5) node {2};
	   \filldraw (\x+10,0) -- (\x+11,1) -- (\x+11,0) -- cycle ;
    }
    \filldraw (18,0) -- (19,1) -- (19,0) -- cycle ;
    \draw (18.75,0.35) node {1};
    \foreach \x in {0,2,4,6}{
	   \draw (\x+10.25,0.65) node {0};
	   \draw (\x+10.75,0.35) node {1};
	   \draw (\x+10.25,2.65) node {3};
        \draw (\x+10.75,2.35) node {4};
	   \draw (\x+10.25,4.65) node {0};
	   \draw (\x+10.75,4.35) node {1};
    }
    \foreach \x in {1,3,5}{
	   \draw (\x+10.25,0.65) node {1};
	   \draw (\x+10.75,0.35) node {0};
	   \draw (\x+10.25,2.65) node {4};
        \draw (\x+10.75,2.35) node {3};
	   \draw (\x+10.25,4.65) node {1};
	   \draw (\x+10.75,4.35) node {0};
    }
	\draw (17.25,0.65) node {1};
    \draw (17.75,0.35) node {0};
    \draw (17.25,2.65) node {4};
    \draw (17.75,2.35) node {3};
    %\draw (\x+0.25,4.65) node {1};
    \draw (17.25,4.65) node {1};

    \foreach \x in {0,2}{
	   \draw (\x+10.25,6.65) node {3};
	   \draw (\x+10.75,6.35) node {4};
        %\draw (\x+0.25,4.65) node {0};
	   \draw (\x+10.75,8.35) node {1};
	   \draw (\x+10,6) -- (\x+11,7);
    }
    \draw (11.25,6.65) node {4};
    \draw (11.75,6.35) node {3};
    \draw (11.75,8.35) node {0};
    \draw (13.25,6.65) node {4};
    \draw (13.75,6.35) node {3};
    \draw (11,6) -- (12,7);
    \draw (13,6) -- (14,7);
    \draw (14.25,6.65) node {3};
    \draw (14.75,6.35) node {4};
    \draw (14,6) -- (15,7);
    \foreach \x in {0,1,2,3}{
	   %\draw (\x,0) -- (\x+1,1);
	   \draw (\x+10.5,5.5) node {2};
	   \draw (\x+10.5,7.5) node {2};
    }
    \draw (14.5,5.5) node {2};
    \draw (15.5,5.5) node {2};
    \draw (19.8,0.5) node {\dots};
    
    \end{tikzpicture}
    
    \caption{The pattern of type $D_r$ (left) and a Young wall in the pattern of type~$D_4$ (right)}
    \label{fig:pattern-D}
\end{figure}

The \textit{pattern of type $D_r$}, shown in the left-hand side of Figure \ref{fig:pattern-D}, consists of two types of boxes: half-boxes carrying possible labels $i \in \{0,1,r-1, r\}$, and full boxes carrying possible labels $1<i<r-1$.
A \textit{Young wall of type $D_r$} is a subset $W$ of the pattern of type $D_r$, satisfying the following rules.

\begin{enumerate}
    \item[(YW1)] $W$ contains all grey half-boxes, and a finite number of the white boxes and half-boxes. 
    \item[(YW2)] $W$ consists of continuous columns of boxes, with no box placed on top of a missing box or half-box. 
    \item[(YW3)] Except for the leftmost column, there are no free positions to the left of any box or half-box. Here the rows of half-boxes are thought of as two parallel rows; only half-boxes of the same orientation have to be present.
    \item[(YW4)] A full column is a column with a full box or both half-boxes present at its top; then no two full columns have the same height.\footnote{Condition (YW4) is called properness, and the arrangements satisfying (YW1)--(YW4) {\em proper Young walls}, in~\cite{kang2004crystal,kwon2006affine}. As we do not consider non-proper Young walls, we omit the adjective {\em proper} for brevity.}
\end{enumerate}

%Scrapped example

Let $\DPart[r]$ denote the set of all Young walls of type $D_r$. In the same way as in the type $A$ case, for $W \in \DPart[r]$ we denote by $\wt_i(W) \in \N$ the number of white boxes or half-boxes of label $i$. The multi-variable generating series of the set $\DPart[r]$ is $$Z_{D_r}(q_0, \ldots, q_r) = \sum_{W \in \DPart[r]} \prod_{i=0}^r q_i^{\wt_i(W)} \; .$$

Define a {\em bar} to be a connected set of boxes and half-boxes, with each half-box occurring once and each box occurring twice. A Young wall $W \in \DPart[r]$ will be called a {\em core} Young wall, if no bar can be removed from it without 
violating the Young wall rules. For an example of bar removal, see \cite[Example 5.1(2)]{kang2004crystal}. Let $\DCorePart{r} \subset \DPart[r]$ denote the set of all 
core Young walls. The following result is the analogue for type $D$ of the Littlewood decomposition, Theorem~\ref{thm:littlewood}. 
\begin{comment}
Removal of a bar can be traced on the abacus similarly to the type A case. However, in the type D there are differences. In fact, the following steps are allowed. For rulers $R_k$, $k \not\equiv 0\; \textrm{ mod}\; (r-1)$:
\begin{enumerate}
 \item[(B1)] If $b$ is a bead at position $s>2r-2$, and there is no bead at position $(s-2r+2)$, then move $b$ 
one position up and switch the color of the beads at 
positions $k$ with $k \equiv 0\; \textrm{ mod}\; (r-1)$, $s-2r+2 < k < s$.
 \item[(B2)] If $b$ and $b'$ are beads at position $s$ and $2r-2-s$ ($1 \leq s \leq r-2$) respectively, then 
remove $b$ and $b'$ and switch the color of the beads 
at positions $k \equiv 0\; \textrm{ mod}\; (r-1), s< k < 2r-2-s$.
\end{enumerate}
For rulers $R_k$, $k \equiv 0\; \textrm{ mod}\; (r-1)$:
\begin{enumerate}
\item[(B3)] Let $b$ be a bead at position $s\geq 2r-2$. If there is no bead at position $(s-r+1)$, and the beads at position $(s-2r+2)$ are of the same color as $b$, then shift $b$ up to position $(s-2r+2)$.
 \item[(B4)] If $b$ and $b'$ are beads at position $s\geq r-1$, then move them up to position $(s-r+1)$. If $s-r+1>0$ and this position already contains beads, then $b$ and $b'$ take that same color.
\end{enumerate}
\end{comment}
%Based on the calculations of \cite{kang2004crystal, kwon2006affine} the following result was obtained in \cite[Proposition 7.2]{gyenge2018euler}. For completeness we include also its proof.
\begin{theorem}[{\cite[Proposition~7.2]{gyenge2018euler}, see also \cite[Proposition~3.2]{gyenge2021young}}]

\label{thm:dncoredecomp}
Given a Young wall 
$W \in \DPart[r]$, any complete sequence of bar removals through Young walls results in the same core $\mathrm{core}(W)\in \DCorePart{r}$, defining a map of sets
\[ \mathrm{core} \colon \DPart[r] \to \DCorePart{r}.\]
%The process can be described on the abacus, respects the decomposition~\eqref{decomp_D_YW}, and results
%in 
There is a combinatorial bijection
\begin{equation*} \DPart[r] \longleftrightarrow  {\mathcal P}^{r+1}  \times \DCorePart{r} \label{Dpart_biject}\end{equation*}
that is compatible with the map $\mathrm{core}$, where ${\mathcal P}$ is the set of ordinary partitions.
\end{theorem}
Let
$$Z_{D_r, \text{core}}(q_0, \ldots, q_r) = \sum_{W \in \DCorePart{r}} \prod_{i=0}^r q_i^{\wt_i(W)}$$
be the generating function of core Young walls. Similarly to the type $A$ case, one can give an explicit expression for the generating function of Young walls using the core-quotient decomposition, and and additional bijection. We will use the change of variables
\[q = q_0q_1q_2^2\cdots q_{r-2}^2q_{r-1}q_r,\]
which once again accounts for the box content of a single bar. 
\begin{theorem}[Kang and Kwon {\cite{kang2004crystal}}]\label{thm:d:fn}\ 
\begin{enumerate}\item 
There is a bijection
\begin{equation*}\label{Dcore_biject} \DCorePart{r} \longleftrightarrow \Z^r.\end{equation*}
\item 
    The generating function of core Young walls of type $D_r$ is given by $$Z_{D_r, \text{core}}(q_0, \ldots, q_r)  = \sum_{w \in \Z^r} q^{\frac{1}{2} w^{\top} C_{D_r} w}\prod_{i=1}^r q_i^{w_i}\; ,$$ where $$C_{D_r} = \begin{pmatrix} 2 & -1 & 0 & & \\ -1 & 2 & \ddots & 0 & 0 \\ 0 & \ddots & \ddots & -1 & -1 \\ & 0 & -1 & 2 & 0  \\ & 0 & -1 & 0 & 2 \end{pmatrix}$$ is the $r\times r$ Cartan matrix of {\em finite} type $D_r$.
\item   The generating function of all Young walls of type $D_r$ is given by $$Z_{D_r}(q_0, \ldots, q_r) = (Z_0)^{r+1}(q) Z_{r, \text{core}}(q_0,\ldots, q_r) = \frac{\displaystyle\sum_{w \in \Z^r} q^{\frac{1}{2} w^{\top} C_{D_r} w}\prod_{i=1}^r q_i^{w_i}}{\displaystyle\prod_{k = 0}^{\infty} (1-q^k)^{r+1}} \; .$$
    \end{enumerate}
\end{theorem}

\subsection{Stratifying the equivariant Hilbert scheme in Type D}
 
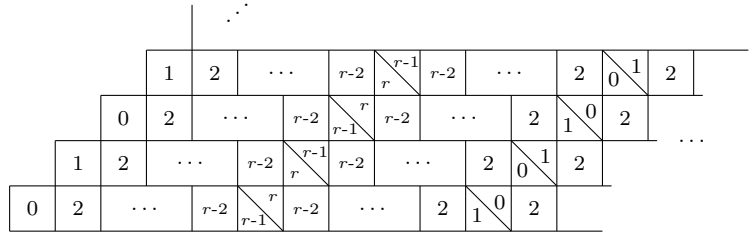
\begin{figure}[h]
\begin{tikzpicture}[scale=0.6, font=\footnotesize, fill=black!20]
  \draw (5,5) node {\reflectbox{$\ddots$}};
  \draw (15,2) node {\dots};
  \draw (4,4) -- (4,5);
  \draw (0,0 ) -- (13,0);
    \foreach \y in {1,2,3,4}
        {
          \draw (\y -1,\y ) -- (\y +12+0.2,\y);
        }
    \foreach \y in {0,1,2,3}
    {
    \foreach \x in {0,1,2,4,5,6,7,9,10,11,12}
                  {
                       \draw (\y+\x ,\y) -- (\y+\x ,\y+1);
                  }
    	%\draw (\x,0) -- (\x+1,1);
    	\draw (\y+6,\y) -- (\y+5,\y+1);
    	\draw (\y+11,\y) -- (\y+10,\y+1);
    	\draw (\y+1.5,\y+0.5) node {2};
    	\draw (\y+4.5,\y+0.5) node {{\tiny $r$-$2$}};
    	\draw (\y+6.5,\y+0.5) node {{\tiny $r$-$2$}};
    	\draw (\y+9.5,\y+0.5) node {2};
    	\draw (\y+11.5,\y+0.5) node {2};
    	\draw(\y+3,\y+0.5) node {\dots};
    	\draw(\y+8,\y+0.5) node {\dots};
    	%\filldraw (\x,0) -- (\x+1,1) -- (\x+1,0) -- cycle ;
    }
     \foreach \y in {0,2}
        {
      	    \draw (\y+0.5,\y+0.5) node {0};
       	    %\draw (\x+0.75,0.35) node {1};
        	\draw (\y+5.36,\y+0.25) node  {{\tiny $r$-$1$}};
        	%\draw (\x+0.49,5.28) node %[fill=white] {{\tiny $r$-$1$}};
        	\draw (\y+5.78,\y+0.75) node {{\tiny $r$}};
        	\draw (\y+10.75,\y+0.65) node {0};
        	\draw (\y+10.25,\y+0.35) node {1};
        }
        \foreach \y in {1,3}
          {
              	    \draw (\y+0.5,\y+0.5) node {1};
               	    %\draw (\x+0.75,0.35) node {1};
                	\draw (\y+5.25,\y+0.25) node  {{\tiny $r$}};
                	%\draw (\x+0.49,5.28) node %[fill=white] {{\tiny $r$-$1$}};
                	\draw (\y+5.70,\y+0.75) node {{\tiny $r$-$1$}};
                	\draw (\y+10.75,\y+0.65) node {1};
                	\draw (\y+10.25,\y+0.35) node {0};
          }     
\end{tikzpicture}
\caption{The transformed pattern of type $D_r$}
\label{fig:transformed:Young}
\end{figure}

Using a very simple transformation, the pattern of type $D_r$, introduced in the previous section, can be related to the 
$G_{D_r}$-equivariant geometry of the affine plane~$\Aff^2$. Following~\cite{gyenge2018euler}, define the
\textit{transformed pattern of type~$D_r$} to be the one presented on Figure~\ref{fig:transformed:Young}; this is related to the pattern of type $D_r$ by an obvious linear transformation.
%The transformation is an affine one, involving a shear: reflect the original Young wall pattern in the line $x = y$ in the plane, translate the $n$th row by $n$ to the right, and remove the grey triangles of the original pattern. 
As it can be checked readily, this pattern is a representation of the coordinate algebra $R=\C[\Aff^2]=\C[x,y]$ of the affine plane and its decomposition into $G_{D_n}$-representations. For $1<i<r-1$, a full box
labelled $j$ below the diagonal, together with its mirror image,
correspond to a 2-dimensional vector space isomorphic to $V_j$ in representation $\rho_j$; a basic example of this is $V=\langle x,y\rangle$ corresponding to representation $\rho_2$. For $j\in\{0,1,r-1,
r\}$, a full block labelled $j$ on the diagonal, as well as a half-block
labelled $j$ below the diagonal with its mirror
image, correspond to a one-dimensional representation.  For example, $\langle 1\rangle \cong \langle x^2y^2\rangle\cong V_0$ are invariants, $\langle xy\rangle \cong \langle x^3y^3\rangle\cong V_1$ correspond to the sign representation, whereas %\footnote{\balazs Is this correct? maybe there are some $i$s in here?}
\[ \langle x^{r-2}, y^{r-2}\rangle \cong \langle x^{r-2}+i^ry^{r-2}\rangle \oplus \langle x^{r-2}-i^ry^{r-2}\rangle \cong \rho_{r-1}\oplus \rho_r. 
\]
The following is the basic result that replaces Theorem~\ref{thm:typeAmain}(2), and connects the combinatorics of type $D_r$ Young walls to the equivariant Hilbert scheme. 
\begin{theorem}[{\cite[Theorem~4.1]{gyenge2018euler}}]
\label{thm:dnorbcells}
For the subgroup $G_{D_r}< SL(2,\C)$ of type $D_r$, 
there is a locally closed decomposition
\[\Hilb(\Aff^2)^{G_{D_r}} = \displaystyle\bigsqcup_{W \in{\mathcal W}_r } \Hilb(\Aff^2)^{G_{D_r}}_W\]
of the equivariant Hilbert scheme $\Hilb(\Aff^2)^{G_{D_r}}$ into strata indexed bijectively by the
set ${\mathcal W}_r$ of Young walls of type $D_r$, with each stratum $\Hilb(\Aff^2)^{G_{D_r}}_W$ isomorphic to a non-empty affine space. 
\end{theorem}
Finding the stratification of $\Hilb(\Aff^2)^{G_{D_r}}$ is relatively straighforward~\cite[Sections 3.5, 4.1]{gyenge2018euler}; the strata $\Hilb(\Aff^2)^{G_{D_r}}_W$ can be thought of as generalised Schubert cells, defined essentially in terms of linear algebra. To describe the geometry of the strata is much harder~\cite[Sections 4.2-4.5]{gyenge2018euler}.

Using standard facts about topological Euler characteristics, this result implies the type $D$ analogue of Theorem~\ref{thm:typeAmain}(3).

\begin{corollary} There is an equality of generating functions 
\[\sum_{(n_i)\in \N^I} \chi_\tp\left(\Hilb^{\oplus_{i\in I} \rho_i^{n_i}}(\Aff^2)\right) \prod_{i\in I}q_i^{n_i} = Z_{D_r}(q_0, \ldots, q_r),
\]
where on the left hand side we consider components of the equivariant Hilbert scheme $\Hilb(\Aff^2)^{G_{D_r}}$.
\end{corollary}

A bijective characterisation of components of the equivariant Hilbert scheme in terms of root system data was given recently in~\cite[Thm.~2]{Paegelow}. 

\subsection{Quot schemes and the substitution formula in Type D}

Given a nonempty subset $J\subset I$, recall the Quot scheme $\Quot_{G_{D_r},J}(\Aff^2)$ of Definition-Theorem~\ref{defthm:quotscheme}, and the corresponding generating function
$$Z_{D_r,J}(q_i\colon i\in J) = \sum_{v_J \in \N^J} \chi_\tp(\Quot_{G_{D_r}, J}^{v_J}(\Aff^2)) \prod_{j \in J} q_j^{v_j}.$$ 
The generating series $Z_{D_r,J}$ can again be obtained via substitution from the Hilbert scheme generating function $Z_{D_r}$. The type $D$ substitution formula is in some sense analogous to the type $A$ formula of Theorem \ref{thm:subst-J-young-walls}, but we have to take a few more things into account. First, a connected component of the subdiagram supported on $I \setminus J$ may now be of finite type $D$. Instead of the denominators $r_i+2$, we thus use the more general expression $h_i+1$, where $h_i$ is the \textit{dual Coxeter number}\footnote{Recall that the (dual) Coxeter number takes the value $(r+1)$ in type~$A_r$, and the value $(2r-2)$ in type~$D_r$.} of the finite type Dynkin diagram in question. Second, even when the connected component is of type $A$, it might connect to the rest of the diagram in vertices other than the endpoints. The necessary information is encoded in integers $c_i$ for each $i \in I \setminus J$, depending on the position of vertex $i$ relative to its finite type component, defined as follows: for $i \in I \setminus J$, let $c_i$ be the sum of the entries in the $i$-th row of the inverse of the finite type Cartan matrix. Third, any vertex in $J$ could now be adjacent to up to three vertices in $I \setminus J$. The substitution formula then reads as follows.
%But in the type D case, the substitution is a bit more complex than in Theorem~\ref{thm:subst-J-young-walls}, because the valency of a node can be 3. {\adam For $i \in I \setminus J$, let $h_i$ be the dual Coxeter number of the connected component of $i$ in the diagram of $I \setminus J$ (note that each such connected is a finite type Dynkin diagram). 
%(the multiplicity of a node is 1 if $i \in \{0,1,n-1,n\}$ and 2 otherwise). 
%For a node $i \in J$, let $k_i, l_i, m_i$ be the value of $h_j$ of the three nodes adjacent to $i$ (or zero if there are less than three nodes connected to $i$). Let moreover $c_i,d_i,e_i$  be the sums of the rows corresponding to the above three adjacent nodes in the inverse Cartan matrix of the respective finite type Dynkin diagrams (or, again, zero if there are less than three nodes connected to $i$).}

\begin{theorem}[\cite{BGS}] \label{thm:subst-J-young-wallsD}
    Given a non-empty subset $J \subset I$, consider the following substitution of variables $q_I$:
    %{\adam \begin{equation}\label{SubstD}
     %   q_i \mapsto s_{r, J}(q_i) = \begin{cases} q_i \exp \left(c_i\frac{2 \pi \sqrt{-1}}{k_i + 1} +d_i \frac{2 \pi \sqrt{-1}}{l_i + 1} +e_i \frac{2 \pi \sqrt{-1}}{m_i + 1} \right), & i \in J, \\ \exp \left(\frac{2 \pi \sqrt{-1}}{h_i + 1} \right), & i \notin J. \end{cases}
    %\end{equation}}
    \begin{equation}\label{SubstD}
        q_i \mapsto s_{r, J}(q_i) = \begin{cases} q_i \exp \left(- \sum_{i \to j \notin J} c_j \frac{2 \pi \sqrt{-1}}{h_j+1}\right), & i \in J, \\ \exp \left(\frac{2 \pi \sqrt{-1}}{h_i + 1} \right), & i \notin J, \end{cases}
    \end{equation}
    where the sum is over arrows in the double quiver from a given $i \in J$ to any $j\in I \setminus J$. Then there is another complex root of unity $c_{r,J}$ such that 
    $$Z_{D_r,J}(q_i\colon i\in J)= c_{r,J} s_{r,J}(Z_r(q_0,\ldots, q_r)).$$
    %Let $m$ denote the minimum of $I$ in $\{0, \ldots, n\}$. Then $$c_{n,I} = (-1)^m \exp \left( 2 \pi \sqrt{-1} \frac{m(m+1)}{2(l_m + 2)}\right)$$
\end{theorem}
In the special case $J=\{0\}$, corresponding to the Hilbert scheme of points of the singular surface $X=\Aff^2/G_{D_r}$ (see Example~\ref{ex:quots}(2)), this result was proved in~\cite{gyenge2018euler}, relying on explicit geometric arguments and the combinatorics of a pattern of type $D_{r,\{0\}}$ (see~\cite[Section~7.3]{gyenge2018euler} in particular).  For different subsets $J\subset I$, this approach appears to be very cumbersome, and the details have not been worked out. The proof in~\cite{BGS} uses instead the method of Nakajima~\cite{NakajimaProof}, which we explain in the next section. 

\subsection{Type E}

Beyond types~$A$ and~$D$, Theorem~\ref{thm:classification} allows only three further subgroups $G<\SL(2,\C)$, the binary exceptional groups of types~$E_6$, $E_7$ and~$E_8$. As in the type $D$ case, only a one-dimensional
torus $T_0\cong \mathbb{C}^\ast$ commutes with the action of the finite group $G_{E_r}$ on $\Aff^2$, with 
non-isolated fixed loci in the Hilbert scheme. On the other hand, the appropriate generalisation of the right hand sides of the formulas in Corollary~\ref{cor_Zr} and Theorem~\ref{thm:d:fn}(3) make sense, defining a series $Z_{E_r}(q_0,\ldots, q_r)$ purely in terms of the Cartan matrix. In~\cite{gyenge2018euler}, a conjecture was formulated which expressed the generating function of Euler characteristics of the Hilbert scheme $\Hilb(\Aff^2/G_{E_r})$, in other words the Quot schemes $\Quot_{G_{E_r},{0}}(\Aff^2)$, in terms of a root-of-unity substitution into this series. Because no technology analogous to type~$D$ Young walls was available, the combinatorial approach did not extend to this case. 

To circumvent this, Nakajima~\cite{NakajimaProof} took a different approach, and analyzed the fibers of the degeneration map~\eqref{eq:quot-degenerationspec}
\[p_{I,\{0\}}:\Hilb(\mathbb{A}^2)^{G} \to \Hilb(\Aff^2/G)\]
from Proposition~\ref{prop:deg}. Using methods relying on the finite-dimensional representation theory of quantum affine Lie algebras and in particular the notion of {\em quantum dimension}, beyond the scope of this review, he was able to give a proof of the substitution formula of~\cite{gyenge2018euler} in all types. Building on this work, Theorem-\ref{thm:subst-J-young-wallsD} was 
shown to hold\footnote{The general substitution formula may also require the (dual) Coxeter numbers of the root systems $E_6, E_7, E_8$ which are 12, 18 and 30 respectively.}
 in all types in~\cite{BGS}, and so in particular in type~$E$.

In a very recent development, Young walls parametrizing a crystal basis of the appropriate representation of the quantum affine algebra of type $E_r$ for $r=6,7,8$, analogous to the type $A$ and $D$ constructions discussed earlier, were constructed by Laurie in~\cite{laurie2025young}. It is natural to wish to connect Laurie's Young walls to the action of the group $G_{E_r}$ on the ring $R=\C[x,y]$ , in the spirit of Figure~\ref{fig:transformed:Young} and the surrounding discussion. One could further hope that these combinatorial objects parameterise a cell decomposition of the equivariant Hilbert scheme $\Hilb(\mathbb{A}^2)^{G_{E_r}}$. All this remains speculation at present. 

\bibliographystyle{amsplain}
\bibliography{main}

\end{document}